\documentclass[12pt]{amsart}
\usepackage[utf8]{inputenc}
\usepackage{bbm}
\usepackage{amssymb}
\usepackage{comment}
\newcommand{\ti}{\text{i}}
\newcommand{\tr}{\text{tr}}

\usepackage[usenames,dvipsnames]{color}
\usepackage{amsmath,amsthm,amsfonts,amssymb,multicol,amscd,amsbsy}
\usepackage{graphicx}
\usepackage{booktabs}
\usepackage{caption}
\usepackage{array}
\usepackage{enumerate}
\usepackage{url}
\usepackage[english]{babel}
\usepackage{mathtools}
\usepackage[toc,page]{appendix}
\usepackage{verbatim} 
\usepackage{enumitem}
\usepackage{bbm} 
\usepackage[normalem]{ulem} 
\usepackage{bm}
\usepackage{hyperref}
\usepackage{tabularx}
\usepackage{latexsym}
\usepackage[margin=1.15in]{geometry}

\newcommand{\floor}[1]{\left\lfloor #1 \right\rfloor}

\newcommand{\be}{\begin}
\newcommand{\e}{\end}
\newcommand{\beq}{\begin{equation}}
\newcommand{\eeq}{\end{equation}}

\newcommand{\ul}{\underline}
\renewcommand{\l}{\left}
\renewcommand{\r}{\right}

\renewcommand{\d}{\mathrm{d}} 

\renewcommand{\Re}{\mathrm{Re}}
\renewcommand{\Im}{\mathrm{Im}}

\newcommand{\set}[1]{\mathbb{#1}}
\newcommand{\curly}[1]{\mathcal{#1}}

\newcommand{\setof}[2]{\left\{ #1\; : \;#2 \right\}}

\newcommand{\R}{\set{R}}
\newcommand{\C}{\set{C}}
\newcommand{\Z}{\set{Z}}

\newcommand{\om}{\omega}
\newcommand{\Om}{\Omega}
\newcommand{\eps}{\epsilon}

\newcommand{\lam}{\lambda}

\newcommand{\gam}{\gamma}

\newcommand{\sig}{\sigma}
\newcommand{\al}{\alpha}
\newcommand{\de}{\delta}

\newtheorem{thm}{Theorem}[section]
\newtheorem{lm}[thm]{Lemma}

\newtheorem{prop}[thm]{Proposition}

\theoremstyle{definition}
\newtheorem{defn}[thm]{Definition}
\newtheorem*{defn*}{Definition}

\numberwithin{equation}{section}

\theoremstyle{remark}

\def\dotuline{\bgroup
  \ifdim\ULdepth=\maxdimen  
   \settodepth\ULdepth{(j}\advance\ULdepth.4pt\fi
  \markoverwith{\begingroup
  \advance\ULdepth0.08ex
  \lower\ULdepth\hbox{\kern.15em .\kern.1em}%
  \endgroup}\ULon}

\def\dashuline{\bgroup
  \ifdim\ULdepth=\maxdimen  
   \settodepth\ULdepth{(j}\advance\ULdepth.4pt\fi
  \markoverwith{\kern.15em
  \vtop{\kern\ULdepth \hrule width .3em}%
  \kern.15em}\ULon}
\allowdisplaybreaks

\title[Local Law from Diophantine Equations]{A Local Law for Singular Values\\ from Diophantine Equations}
\date{May 8, 2020}
\author[A.\ Adhikari]{Arka Adhikari}
\email{adhikari@math.harvard.edu}
\author[M.\ Lemm]{Marius Lemm}
\email{mlemm@math.harvard.edu}
\address{Department of Mathematics, Harvard University, 1 Oxford Street, Cambridge, MA 02138, USA}

\begin{document}

\begin{abstract}
We introduce the $N\times N$ random matrices
$$X_{j,k}=\exp\left(2\pi i \sum_{q=1}^d\ \omega_{j,q} k^q\right)
\quad \textnormal{with } \{\omega_{j,q}\}_{\substack{1\leq j\leq N\\ 1\leq q\leq d}} \textnormal{ i.i.d.\ random variables},
$$
and $d$ a fixed integer. We prove that the distribution of their singular values converges to the local Marchenko-Pastur law at scales $N^{-\theta_d}$ for an explicit, small $\theta_d>0$, as long as $d\geq 18$. To our knowledge, this is the first instance of a random matrix ensemble that is explicitly defined in terms of only $O(N)$ random variables exhibiting a universal local spectral law. Our main technical contribution is to derive concentration bounds for the Stieltjes transform that simultaneously take into account stochastic and oscillatory cancellations. Important ingredients in our proof are strong estimates on the number of solutions to Diophantine equations (in the form of Vinogradov's main conjecture recently proved by Bourgain-Demeter-Guth) and a pigeonhole argument that combines the Ward identity with an algebraic uniqueness condition for Diophantine equations derived from the Newton-Girard identities.
\end{abstract}

\maketitle

\section{Introduction}
In 1955, Eugene Wigner introduced random matrices drawn from what are now called the Gaussian  Unitary Ensemble (GUE) and Gaussian Orthogonal Ensemble (GOE) as toy models of the \textit{deterministic} quantum Hamiltonians describing heavy nuclei \cite{Wigner}. Wigner noticed that the eigenvalue spacing of the random matrices matched experimental data for the spacing distribution of the energy levels to surprising accuracy. This discovery was subsequently broadened into the highly influential Wigner-Dyson-Mehta-Gaudin universality conjecture which says that the eigenvalue spacing distribution of a matrix ensemble depends only on its symmetry class. In the past 15 years, a number of celebrated results in random matrix theory succeeded in verifying the Wigner-Dyson-Mehta-Gaudin conjecture in great generality \cite{ERetal,ESY2,ESYY,TV}.

An avenue of investigation which still remains to be fully understood is the fact that the universality of random matrix statistics appears to extend to various deterministic systems. The first observation of this kind was made in Wigner's foundational study \cite{Wigner} of heavy nuclei. Two other famous examples are Montgomery's Pair Correlation Conjecture (zeros of the Riemann zeta function follow GUE statistics) and the Quantum Chaos Conjecture (Laplace eigenvalues on classically chaotic domains follow GUE statistics). More generally, random matrix statistics constitute an extremely wide-ranging universality class for highly correlated point processes which empirically appears to include a number of deterministic and real-world examples. Therefore, it is a central goal of modern research in the field to derive random matrix eigenvalue statistics for matrix ensembles with as little randomness as possible. 

In the present paper, we are not interested in universality of eigenvalue spacing (which is a very fine statement that requires understanding individual eigenvalues), but instead we study the convergence of the empirical spectral distribution $\frac{1}{N}\sum_{j=1}^N \de_{\lam_j}$, where $\{\lam_j\}_{1\leq j\leq N}$ are the eigenvalues of the matrix under investigation. Understanding its behavior down to small scales is a fundamental ingredient to all proofs of universality of eigenvalue spacing. More precisely, one aims to prove the weak convergence of $\frac{1}{N}\sum_{j=1}^N \de_{\lam_j}$ to a well-defined (and also appreciably universal) limiting distribution, a statement that can be seen as a non-commutative analog of the central limit theorem. The most famous limiting distribution is the Wigner semicircle law which arises for general ensembles of Hermitian matrices \cite{Wigner}. When the weak convergence is proved with respect to order-one test functions, such a statement is called a \textit{global law} for the empirical spectral distribution. A refinement where the test functions live on scales $N^{-\theta}$ (so scales shrinking with $N$) is instead called a \textit{local law}. One can only expect this for $\theta<1$ because the typical eigenvalue spacing is $N^{-1}$. For further background on local laws, we refer to the books \cite{BK,EY}.

In our recent work \cite{ALY}, we considered a novel ensemble of random matrices which is rather structured: all entries in a given row are fully dependent. Specifically, each row is obtained by evaluating the complex exponential along orbits of the skew-shift $\binom{j}{2}\omega+jy+x \text{ mod } 1$ with $\om$ an irrational parameter. The main result of \cite{ALY} establishes a \textit{global} law when the initial values of $y$ (a starting coordinate of the skew-shift) for every row are i.i.d.\ uniform random variables. The basic idea is that the oscillations coming from the complex exponentials with irrational frequency end up supporting the comparatively small amount of randomness. In this way, \cite{ALY}  establishes a global law for a random matrix ensemble comprised of ``only'' $N$ independent random variables. This count of $N$ is to be compared to the classical ensembles of random matrix theory which hold order $N^2$ independent random variables.

It was left as an open problem in \cite{ALY} to derive the first \textit{local} law for a random matrix ensemble in which all entries in a given row are fully dependent and which thus depends on only order $N$ random variables. This open problem is addressed in the present paper (Theorem \ref{thm:main1}). We view our result as a step forward in the important long-term program of deriving random matrix statistics for systems that are progressively less random and more structured. The method also yields weak delocalization bounds for eigenvectors (Theorem \ref{thm:deloc}), another hallmark of random matrix behavior. The matrix ensemble \eqref{eq:XNdefn} we propose here is inspired by the skew-shift ensemble from \cite{ALY} but incorporates higher-degree polynomial terms. As in \cite{ALY}, the main technical challenge is to harness stochastic and oscillatory cancellations hand-in-hand. This involves combining techniques from probability theory, harmonic analysis, and number theory, with the latter arising from the close connection between resonances of exponential sums and Diophantine equations. We expect that the ensemble which we propose here has even better properties, namely a local law all the way down to the nearly optimal scale $N^{-1+\eps}$ and universal gap statistics (see Conjecture \ref{conj}).

We would like to mention that in recent years numerous works have established local laws have for matrix ensembles with more correlation than Wigner-type matrices. For example, we mention the works on matrices with polynomially decaying correlation structure \cite{EKS} or with correlated random variables generated from statistical physics \cite{FL,HKW,KK}, on adjacency matrices of Erd{\H o}s-Renyi graphs with probability $p\geq N^{-1+\eps}$ \cite{HL} (i.e., random graphs of average degree $pN=N^\eps$) and on adjacency matrices of random regular graphs of large, but fixed degree \cite{BKY}. 

Still, to our knowledge there exists no proof in the vast and constantly growing random matrix theory literature of even a coarse local law for a random matrix ensemble that is explicitly defined in terms of $O(N)$ random variables, a result we supply here. In this regard, it should be mentioned that specifically the adjacency matrices of random regular graphs of fixed degree for which local law and universality were proved in \cite{BKY,BHKY} are also highly structured and arguably hold a comparably small amount of randomness. 

At any rate, regardless of how one precisely quantifies the degree of randomness of these correlated random matrix models, the mechanism underlying the emergence of a local law is novel in the present work and it shows that oscillatory cancellations can effectively mimic stochastic cancellations in the spectral theory of random matrices. 



\section{Model and main results}
\subsection{The model}
Inspired by the polynomial nature of the skew-shift, we consider the following random matrix ensemble.

\be{defn}[The random matrix ensemble]
Let $d$ be an integer parameter and for every integer $N\geq 1$, let
    $$
    \{\omega_{j,q}\}_{\substack{1\leq j\leq N\\ 1\leq q\leq d}}
    $$
    
    be a family of independent, identically distributed random variables on the interval $[0,1]$. Define the $N\times N$ matrix $X_N$ by
\begin{equation}\label{eq:XNdefn}
    [X_N]_{jk} = \frac{1}{\sqrt{N}}\exp\l(2\pi \ti \sum_{q=1}^d  \omega_{j,q}  k^q\r),
\end{equation}
\end{defn}

We note that $X_N$ is indeed constructed from a total of $dN=O(N)$ independent random variables.

 The main result concerns the distribution of the singular values of $X_N$, call them $\sigma_1,\ldots,\sigma_N$ on a local scale (i.e., a scale of the form $N^{-\theta}$). The local law is most conveniently  as the convergence of the \textit{Stieltjes transform} of the empirical spectral measure,
\beq\label{eq:mNdefn}
m_N(z)=\frac{1}{N}\tr\l(\frac{1}{X_N(X_N)^*-z}\r)
=\frac{1}{N}\sum_{j=1}^N \frac{1}{x-z}\delta_{\sigma_j^2}(x) , 
\eeq
for $z\in\C$ with $\Im[z]=N^{-\theta}$. The appropriate limiting distribution for singular values of square matrices is the Marchenko-Pastur distribution \cite{MP} with density parameter equal to $1$, i.e.,
\beq\label{eq:MPdefn}
\rho_{\mathrm{MP}}(x)=\frac{1}{2\pi x}\sqrt{x(4-x)}\mathbbm 1_{0\leq x\leq 4}
\eeq
with Stieltjes transform
\beq
m_{\mathrm{MP}}(z)=\int_{\R} \frac{1}{x-z} \rho_{\mathrm{MP}}(x)\d x
\eeq
defined for all $z\in\C$ with positive imaginary part. 

  \be{rmk}
  We recall that Marchenko-Pastur \cite{MP} showed that $\rho_{\mathrm{MP}}(x)$ arises as the distribution of singular values of sample covariance matrices $X^*X$ when $X$ is comprised of independent and identically distributed Gaussian entries. Note also that $\mathrm{d}\mu_{\mathrm{MP}}(x)$ is the push forward measure of the Wigner semicircle law under the map $x\mapsto x^2$, so this result is in accordance with the semicircle law. 
  \e{rmk}
    
    \subsection{Main results}

\begin{defn}\label{defn:whp}
We say that the estimate $A\leq B$ holds \textit{with high probability}, if $\mathbb P(B>A)\to 0$ as $N\to\infty$ for every fixed value of the parameter $d$.
\end{defn}

Our notion of ``with high probability'' is not quantitative and obtaining strong probabilistic estimates is not our focus here. The proof is quantitative and can easily yields explicit polynomial bounds on the relevant small probabilities if desired, but in contrast to the case of Wigner matrices one does not get arbitrarily large polynomial decay of the probability of failure.

We come to our assumptions on the distribution of the $dN$ i.i.d.\ random variables $\{\om_{j,q}\}$. First, we may restrict their support to the $1$-torus $\R/\Z$ ($=[0,1]$ with endpoints identified) without loss of generality because of the complex exponential in the definition of the matrix model \eqref{eq:XNdefn}. We then make the following convenient regularity assumption.

\be{ass}\label{ass:omega}
Let $\rho:\R/\Z\to [0,\infty)$ be a probability density function on the 1-torus with bounded derivative, $\|\rho'\|_\infty\leq C$. For every integer $N\geq 1$, 
    $$
    \{\omega_{j,q}\}_{\substack{1\leq j\leq N\\ 1\leq q\leq d}}
    $$
    is a family of independent random variables on $\R/\Z$ each with distribution $\rho(\om)\d \om$
\e{ass}

The derivative bound from Assumption \ref{ass:omega} is only used in the proof of Lemma \ref{lm:omegaexp} and can be weakened if desired. We also mention in passing that the argument simplifies slightly if the uniform distribution is used throughout, i.e., if $\rho(x)=1$, but the main challenges stay the same.

 We introduce the $d$-dependent variables
\beq\label{eq:thetap0defn}
\theta_0:=\frac{\frac{p}{18}-1}{2p+4},\qquad \text{with } p:=\floor{\frac{17d}{16}}.
\eeq 

Our main result is the following local law, formulated in terms of Stieltjes transforms.

\begin{thm}[Main result]\label{thm:main1}Let $d\geq 18$ and $0<\theta<\theta_0$. 

Then, for every $\kappa\in(0,1)$ there exists $c_\kappa>0$ so that on the domain
\beq
\curly{D}:=\setof{z=E+i\eta\in \C}{E\in (\kappa,4-\kappa) \textnormal{ and } \eta\in (N^{-\theta},c_\kappa) },
\eeq
it holds that
\begin{equation}\label{eq:main}
    \sup_{z\in\curly{D}}|m_N(z) - m_{\mathrm{MP}}(z)| \le \frac{N^{-\theta_0}}{\Im[z]},
\end{equation}
with high probability in the sense of Definition \ref{defn:whp}.
\end{thm}

We make a few remarks about this result.

\be{rmk}
\be{enumerate}[label=(\roman*)]
\item 
The condition $d\geq 18$ ensures that $\theta_0>0$ and thus that the domain $\curly{D}$ is non-empty. 

\item As $d\to\infty$, we have $\theta_0\to \frac{1}{18}$ and so the smallest possible scale on which Theorem \ref{thm:main1} gives a local law is $\Im[z]=N^{-\frac{1}{18}+\eps}$. While it is not our goal to optimize the scale here, we note that improving the local law to the $N^{-1+\eps}$ that holds for most previously studied random matrix ensembles will likely require new ideas.

\item In Theorem \ref{thm:main1}, we made some effort to choose constants that yield reasonable values of $\theta_0$ while applying for all sufficiently large $d$. One may also ask what the smallest value of $d$ is that can be treated by the general method. We found that \eqref{eq:main} holds for $d=10$ (but not for all $d\geq 10$; note that there is inherent non-monotonic behavior coming from taking integer parts), if one modifies the relevant parameters as follows: One takes $p$ to $\floor{\frac{15d}{16}}$, $\theta_0$ to $\frac{\frac{p}{8.99}-1}{2p+4}$ and sets $\theta'$ equal to $0.21$ instead of $\frac{1}{5}$ in the proof of Theorem \ref{thm:concentration}.

\item The techniques extend straightforwardly to energies near the spectral edges, i.e., to $E\approx 0$ and $E\approx 4$, but the lower bound on $d$ will increase. We decided to forgo the details here to keep the focus on the new ideas. 
 
\e{enumerate}
\e{rmk}

\subsection{Rigidity of eigenvalues}
By standard techniques, Theorem \ref{thm:main1} implies rigidity of the singular value distribution in the following way. We recall that $0\le \sig_1^2\le \sig_2^2\le \ldots \le \sig_N^2$ denote the eigenvalues of the matrix $X_NX_N^*$, or equivalently, the squared singular values of $X_N$. We define the cumulative empirical spectral distribution (or eigenvalue counting function, or integrated density of states) by
\begin{equation}
    F_N(E)=  \frac{1}{N} \sum_{j=1}^N \mathbbm{1}(\sig_j^2 \le E).
\end{equation}
The appropriate limiting object is the eigenvalue counting function for the Marchenko-Pastur distribution,
\begin{equation}
    F_{\mathrm{MP}}(E) = \int_{-\infty}^{E} \rho_{\mathrm{MP}}(x) \text{d}x,
\end{equation}
with $\rho_{\mathrm{\mathrm{MP}}}(x)$ defined in \eqref{eq:MPdefn}.

\be{cor}[Rigidity of eigenvalues]\label{cor:rigidity}
For every $E\in (\kappa, 4-\kappa)$ and every $0<\theta<\theta_0$, 
\beq\label{eq:corcounting}
|F(E)-F_{\mathrm{\mathrm{MP}}}(E)|\leq N^{-\theta}
\eeq
holds with high probability in the sense of Definition \ref{defn:whp}.
\e{cor}

The rigidity estimate \eqref{eq:corcounting} clarifies that it is natural to have the factor $\frac{1}{\Im[z]}$ on the right-hand side in \eqref{eq:main}. Corollary \ref{cor:rigidity} is proved in the appendix.

\subsection{Delocalization bounds for eigenvectors}

The proof strategy behind Theorem \ref{thm:main1} is described in the next subsection. It provides detailed information on the diagonal of the Green's function. From these we can readily conclude the following delocalization bounds for eigenvectors with respect to the canonical basis. While these bounds are relatively weak due to the constraint on $\theta_0$ in Theorem \ref{thm:main1}, they highlight another way in which the model \eqref{eq:XNdefn} behaves similarly to matrix ensembles with many more random variables. 

We write $u_1,\ldots,u_N$ for a choice of $N$ linearly independent eigenvectors of $X_N X_N^*$ which are $\ell^2$-normalized, i.e., $\sum_{i=1}^N |u_\al(i)|^2=1$. We measure localization via the $\infty$-norm of each $u_\al$. As benchmarks, we note that in this normalization a sharply localized vector satisfies $\|u_\al\|_\infty \geq c>0$ with $c$ independent of $N$, while at the other end of the scale, a fully delocalized vector satisfies  $\|u_\al\|_\infty \leq cN^{-1/2}$. The bound we prove here says that the eigenvectors are at least weakly delocalized.

\be{thm}[Delocalization bound for eigenvectors]\label{thm:deloc}
Let $H_N u_\al=Eu_\al$ with $E\in (\kappa,4-\kappa)$ for some $\kappa \in (0,1)$. Let $0<\theta<\theta_0$. Then
$$
\max_{1\leq i\leq N} |u_\al(i)|^2\leq C N^{-\theta} 
$$
holds with high probability in the sense of Definition \ref{defn:whp}.
\e{thm}

This result is a ``corollary of the proof'' of Theorem \ref{thm:main1}. The argument is deferred to the appendix.

\subsection{Proof Strategy}
The effective spectral analysis of the ensemble \eqref{eq:XNdefn} relies crucially on oscillatory cancellations, in the form of exponential sum estimates, replacing the role played by stochastic cancellations in the more traditional probabilistic proofs. This was also the case in \cite{ALY}, but the present case is more delicate because one no longer has access to the moment method when proving a local law.

There exists by now a well-established strategy for deriving local laws via the Stieltjes transform. We write $z=E+i\eta$ and note that $m_{\mathrm{MP}}(z)$ is characterized as the unique solution to the quadratic equation
\beq\label{eq:mpquadratic}
zm_{\mathrm{MP}}(z)^2 + zm_{\mathrm{MP}}(z)+1=0
\eeq
with positive imaginary part. The main idea is that the Schur complement implies that $m_N(z)$ satisfies an approximate version of \eqref{eq:mpquadratic}. The steps of the proof are as follows (cf.\ Section 8 in \cite{EY}). 

\begin{enumerate}
  \item Schur complement formula and partial expectations
  \item Self-consistent equation for the Stieltjes transform through interlacement of eigenvalues
 \item Large deviation estimate for the fluctuations in the self-consistent equation
\item Initial-scale estimate for large $\eta$
\item Bootstrap argument based on Lipschitz continuity
\end{enumerate}

We modify steps $(1),(2),(4)$ and $(5)$ in various small ways to comply with our specific matrix ensemble and the comparatively large $\eta$-scale that we are working on, e.g., the operator identity \eqref{eq:B*B} in Step (1) is helpful for matrices with independent rows but strong dependencies within rows. However, the real crux lies with step $(3)$. This is expected since step (3) is the main place where the precise stochastic nature of the matrix ensemble enters. For the standard ensembles, one can rely on the well-developed large deviation principles for quadratic forms of i.i.d.\ random variables \`a la Hanson-Wright and controlling the variance through the Ward identity.

However, in the present ensemble \eqref{eq:XNdefn} the relevant quadratic form is still highly correlated (see Proposition \ref{prop:initial} for its moments). Hence, the standard concentration techniques fail completely in our model. 

Remedying this and deriving concentration of this correlated quadratic form is the main technical contribution of the present work. In doing so, we resolve the following two technical challenges: First, we need to connect oscillatory cancellations and stochastic cancellations. A key observation in this regard is that the number of terms after partial averaging can be identified with solutions to certain Diophantine equations and so the number of terms can be controlled via the breakthrough result of Bourgain-Demeter-Guth confirming the Vinogradov main conjecture \cite{BDG}. 

After invoking this powerful result as an a priori estimate on the number of terms that need to be treated, there still remains the second fundamental technical difficulty of controlling the size of the off-diagonal Green's function entries, i.e., $(X_N^*X_N-z)^{-1}_{ij}$ with $i \ne j$. While for Wigner matrices one can show the optimal estimate $\frac{1}{\sqrt{N \eta}}$ on these entries by an induction, this induction argument breaks down in our (much more dependent) model. On top of this, the types of Green's function entries that are relevant to the Wigner moment computation have a very regular structure due to the stochastic cancellations that occur, while in our case they are labeled by solutions to Diophantine equations. This makes it initially completely unclear (even on a heuristic level!) why the off-diagonal Green's function entries should be sufficiently small, beyond the weak a priori bounds implied by the Ward Identity. 

We address this problem by deriving an algebraic uniqueness condition for the relevant Diophantine equations from the Newton-Girard identities (Lemma \ref{lm:uniq}) and combining this algebraic fact with the Ward identity to devise a rather delicate pigeonhole principle controlling the number of ``large'' Green's function entries. The resulting procedure leverages the relatively small gain from the pigeonhole principle and yields a local law on a fixed, but relatively large scale. 

\subsection{Discussion}
\subsubsection{Ergodic Theory Background}
Let us explain the choice of matrix ensemble \eqref{eq:XNdefn}. For this, it is beneficial to take an ergodic theory perspective. The skew-shift dynamics that was used to generate the matrix model in \cite{ALY} and that directly inspired the present choice \eqref{eq:XNdefn} has gained notoriety in other contexts as an ergodic dynamical system which is able to generate quasi-random behavior despite being almost as rigid as irrational circle rotation. While quasi-random behavior of the skew-shift is not yet fully understood, relevant partial results exist in the study of one-dimensional Schr\"odinger cocycles with potential obtained by sampling the cosine along the skew-shift \cite{B1,B2,BGS,HLS1,HLS2,KL,K1} and the famous Poissonian conjecture of Rudnick-Sarnak-Zaharescu \cite{RSZ} in \cite{DRH,MS,MY,RS}. 

Comparing with the situation for one-dimensional Schr\"odinger operators, the fact that our results here improve with increasing $d$ can be seen as a random-matrix analog of a result of Kr\"uger \cite{K1} that the potential generated from sampling the cosine along the $d$-dimensional skew-shift (with $d$ large) has positive Lyapunov exponent at small coupling. 

If one is interested in a local law for a fully dynamically generated random matrix model, we mention that \eqref{eq:XNdefn} could conceivably be replaced with the orbits of skew-shifts on $d$-dimensional tori (with $d$ large but finite) in which case the random variables could be interpreted as a random starting position of the dynamical system.

\subsubsection{Future Directions}
We emphasize that local laws (more specifically the spectral rigidity that they imply) constitute Step 1 in the well-known 3-step strategy of Erd{\H o}s-Schlein-Yau for proving universality of the gap distribution of various random matrix ensembles \cite{EY}. Hence, the local law proved here is not only itself an instance of random matrix statistics in a not-so-random ensemble, the local law may also have a role to play in proving the universality of the gap distribution of the model \eqref{eq:XNdefn}, which can be observed numerically. Implementing this will first require improving the scale of the local law, i.e., increasing the relatively small $\theta_0$ found in Theorem \ref{thm:main1}.

In view of the considerations above, we propose the following conjecture about the key parameters $\theta_0$ and $d$ in Theorem \ref{thm:main1}.

\be{conj}\label{conj}
For every $d\geq 2$ and $\eps>0$, the local law \eqref{eq:main} holds with $\theta_0=1-\eps$.
\e{conj}

Our belief that $d\geq 2$ should be the correct condition is in line with the widespread heuristic that the skew-shift on the $2$-tous (which also involves a quadratic nonlinearity) displays random-like spectral behavior in other contexts as reviewed above, while the quasiperiodic case $d=1$ does not.

\subsubsection{A possible refinement by iteration}
We close the discussion by sketching an optional method for slightly improving the scale in Theorem \ref{thm:main1} if desired. 
The following kind of bootstrap argument can be used to slightly improve the value of $\theta$ a posteriori. Observe that the proof of Theorem \ref{thm:main1} utilizes the a priori bound on the Green's function that follows from the naive Ward identity (cf.\ Lemma \ref{lm:ap}). Any improvements of this a priori estimate feed forward through the argument and moderately increase $\theta$. One way to obtain an improvement is to use the result itself (more precisely Corollary \ref{cor:rigidity}) in the following fashion. Let $G=(X_N(X_N)^*-z)^{-1}$. By the Ward identity for $\tilde G=(X_N(X_N)^*-z)^{-1}$ and the fact that $X_N(X_N)^*$ and $(X_N)^*X_N$ have the same non-zero eigenvalues, we have
$$
    \sum_{i,j=1}^N |\tilde{G}|^2_{ij} = \frac{\tr [\Im[\tilde{G}]]}{\eta}
    = \frac{\tr [\Im[G]]}{\eta}
    =\sum_{i=1}^N\frac{1}{|\lambda_i -z|^2}.
$$
Since the rigidity established by Corollary \ref{cor:rigidity} bounds the right-hand side, it can be used to obtain an improved pigeonhole principle for $\tilde G$. Implementing row-removal appropriately through eigenvalue interlacement,  one can derive a modest improvement to the scale of the local law. We leave the details to the interested reader.

\section{Step 1: Schur complement formula and partial expectation}
In this section, we adapt the standard derivation of the self-consistent equation via the Schur complement formula and the interlacement of eigenvalues between matrices and their minors. Moreover, a crucial role is played by the operator identity
\beq
\label{eq:B*B} 
A^* \frac{1}{AA^*-z} A= A^*A \frac{1}{A^*A-z}
\eeq
which was already used in the random matrix context in \cite{ESYY}. This operator identity is not used in more recent renditions of the method but it is crucial for the present ensemble which has independent rows but strong correlations within each row.

\subsection{Schur complement formula}
Let $1\leq i\leq N$ and use $X_N^{(i)}$ to designate the $(N-1)\times N$ matrix that equals $X_N$ with the $i$th row removed. We also denote $m_N^{(i)}(z)=\tr[(X_N^{(i)}(X_N^{(i)})^*-z)^{-1}]$ and write $r^i$ for the $i$th row of the matrix $X_N$, with the convention that $r^i$ is represented as a column vector in $\C^N$.
By the Schur Complement formula and $r^i(r^i)^*=1$,
\begin{equation} \label{eq:Schurcomplement}
\begin{aligned}
    [(X_N(X_N)^*-z)^{-1}]_{ii}= \frac{1}{1 - z - r^i (X_N^{(i)})^*(X_N^{(i)}(X_N^{(i)})^*-z)^{-1}X_N^{(i)} r^i }
    \end{aligned}
\end{equation}
We apply the matrix identity \eqref{eq:B*B} with $A=X_N^{(i)}$ to obtain
\begin{equation}
\label{eq:XX*}
    (X_N^{(i)})^*(X_N^{(i)}(X_N^{(i)})^* -z)^{-1}X_N^{(i)} = (X_N^{(i)})^* X_N^{(i)}((X_N^{(i)})^* X_N^{(i)} -z)^{-1}.
\end{equation}
Next we diagonalize the $N\times N$ matrix $(X_N^{(i)})^*X_N^{(i)}$. Note that it has at least one eigenvector with eigenvalue zero since $X_N^{(i)}$ has rank at most $N-1$. We call this eigenvector $e_N^{(i)}$ and denote $\lambda^{(i)}_N=0$. The remaining eigenvalues are $\lambda_1^{(i)} \le \lambda_2^{(i)} \le \ldots \le \lambda_{N-1}^{(i)}$ and the corresponding eigenvectors are $e_1^{(i)},\ldots,e_{N-1}^{(i)}$. Using the eigenbasis, we obtain
\begin{equation}
\begin{aligned}
   (X_N^{(i)})^* X_N^{(i)}((X_N^{(i)})^* X_N^{(i)} -z)^{-1}
    =
    \sum_{j=1}^{N} \frac{\lambda_j^{(i)}|\langle r^i, e_j^{(i)}\rangle|^2}{\lambda_j^{(i)}-z} 
    = 1 + z \sum_{j=1}^{N} \frac{|\langle r^i, e_j^{(i)}\rangle|^2}{\lambda_j^{(i)}-z}.
        \end{aligned}
\end{equation}

\subsection{Partial expectation}
For fixed $i$, we introduce the partial expectation $\mathbb E_i$ with respect to the random variables $\{\omega_{i,q}\}_{1\leq q\leq d}$. Since the $\lambda_j^{(i)}$ and $e_j^{(i)}$ are independent of these random variables and 
$$
r^i(k)=\frac{1}{\sqrt{N}}e\l[\sum_{q=1}^d \om_{i,q} k^q\r],\qquad \textnormal{with } e[\cdot]=\exp(2\pi\ti (\cdot)),
$$
we have
$$
\begin{aligned}
\mathbb E_i \l[\sum_{j=1}^{N} \frac{|\langle r^i, e_j^{(i)}\rangle|^2}{\lambda_j^{(i)}-z}\r]
=&
\sum_{j,k,l=1}^N \frac{e^{(i)}_j(k)e^{(i)}_j(l)}{\lambda_j^{(i)}-z}
\int_0^1 \ldots\int_0^1 \overline{r^i(k)r^i(l)} \d\omega_{i,1}\ldots \d\omega_{i,d}\\
=&\frac{1}{N}\sum_{j=1}^N \frac{1}{\lambda_j^{(i)}-z}=m_N^{(i)}(z)\\
\end{aligned}.
$$
We introduce the fluctuation term
\beq\label{eq:Fidefn}
F_i(z)=
\sum_{j=1}^{N} \frac{|\langle r^i, e_j^{(i)}\rangle|^2}{\lambda_j^{(i)}-z}
-\mathbb E_i \l[\sum_{j=1}^{N} \frac{|\langle r^i, e_j^{(i)}\rangle|^2}{\lambda_j^{(i)}-z}\r]
= \frac{1}{N} \sum_{\substack{k,l=1\\ k \ne l}}^{N} G^{(i)}_{k,l} e\l[ \sum_{q=1}^d(k^q -l^q)\omega_{i,q}\r]
\eeq
where we defined the Green's function
\beq\label{eq:Gidefn}
G^{(i)}=((X_N^{(i)})^{*}X^{(i)}_N -z)^{-1}.
\eeq
Returning to the Schur complement formula \eqref{eq:Schurcomplement}, we have shown that
\beq\label{eq:Schurcomplement'}
 [(X_N(X_N)^*-z)^{-1}]_{ii}= \frac{1}{- z - zm_N^{(i)}(z)-zF_i(z)}.
\eeq

We now sum this identity over $1\leq i\leq N$ and conclude the following result.

\begin{prop}\label{prop:step1}
For $1\leq i\leq N$, define $F_i(z)$ by \eqref{eq:Fidefn}. Then
\beq\label{eq:step1}
m_N(z)+\frac{1}{z}\sum_{i=1}^N \frac{1}{1+m_N^{(i)}(z)+F_i(z)}=0
\eeq
\end{prop}

We can rewrite \eqref{eq:step1} in the form of a self-consistent equation for $m_N(z)$, 
\beq\label{eq:step2}
m_N(z)+\frac{1}{z}\sum_{i=1}^N \frac{1}{1+m_N(z)+\curly{E}_i(z)}=0,\qquad \text{with }\curly{E}_i(z):=m_N^{(i)}(z)-m_N(z)+F_i(z),
\eeq
and $\curly{E}_i(z)$ will be a small error term.

\section{Step 2: Interlacement and the self-consistent equation}
In this short section, we bound the first contribution to $\curly{E}_i(z)$, namely $m_N^{(i)}(z)-m_N(z)$, by a standard argument based on the eigenvalue interlacement of matrices and their minors.

\be{lm}\label{lm:interlace}
There exists a universal constant $C>0$, so that for every $1\leq i\leq N$ and all $z=E+i\eta \in\curly{D}$,
$$
\l|m_N(z)-m_N^{(i)}(z)\r|\leq \frac{C}{N\eta}.
$$
\e{lm}

\be{proof}
We have
$$
m_N^{(i)}(z)=\frac{1}{N}\tr((X_N^{(i)})^{*}X^{(i)}_N -z)^{-1}=\frac{1}{N}\sum_{j=1}^N \frac{1}{\lambda_j^{(i)}-z}
$$
with $\lambda_N^{(i)}=0$. By the min-max characterization of singular values, the non-zero eigenvalues among the $\{\lambda_j^{(i)}\}_{1 \le j \le N-1}$ are also eigenvalues of the $(N-1)\times(N-1)$ matrix $B=X_N^{(i)}(X_N^{(i)})^*$ with the same multiplicity. Moreover, dimension counting implies that $\dim\ker B=\dim\ker((X_N^{(i)})^{*}X^{(i)}_N)-1$ and so
$$
\l|m_N^{(i)}(z)-\frac{1}{N}\tr\l(\frac{1}{B-z}\r)\r|=\frac{1}{N|z|}\leq \frac{1}{N\eta}
$$
Define $H=X_N(X_N)^*$ and observe that $B$ is obtained from $H$ by removing the $i$th row and column. By Cauchy's interlacing theorem,  the eigenvalues $\{\lambda_j^{(i)}\}_{1\leq j\leq N-1}$ and $\{\lambda_j\}_{1\leq j\leq N}$ interlace. From integration by parts, as e.g. in the proof of Lemma 7.5 in \cite{EY}, it follows that 
$$
\l|\frac{1}{N-1}\tr\l(\frac{1}{B-z}\r)-m_N(z)\r|
=\l|\frac{1}{N-1}\tr\l(\frac{1}{B-z}\r)-\frac{1}{N}\tr\l(\frac{1}{H-z}\r)\r|\leq \frac{C}{N\eta}
$$
for a universal constant $C>0$. This proves Lemma \ref{lm:interlace}.
\e{proof}

\section{Step 3: Large deviation bounds}
In this section we control the other contribution to $\curly{E}_i(z)$, the fluctuations $F_i(z)$ defined in \eqref{eq:Fidefn} in the following way. We recall the Definition \eqref{eq:thetap0defn} of $\theta_0$ and $p$.

\begin{thm}[Moment bound]\label{thm:concentration}
Let $d\geq 18$, $\eps>0$ and let $\Im[z]\geq N^{-\theta}$ with $0<\theta<\theta_0$. 
Then
\beq
|\mathbb E_i[F_i(z)^{2p}]|
\leq \l(\frac{N^{-\frac{1}{36}-\eps}}{\Im[z]}\r)^{2p}.
\eeq
holds for sufficiently large $N$.
\end{thm}

Markov's inequality then implies the following large deviation estimate.

\be{cor}[Large deviation estimate]\label{cor:ld}
Let $\beta>0$. Under the assumptions of Theorem \ref{thm:concentration}, 
\beq
\mathbb P\l(|F_i(z)|\geq \frac{N^{\frac{\beta }{2p}-\frac{1}{36}-\eps}}{\Im[z]}\r)\leq N^{-\beta}
\eeq
holds for sufficiently large $N$.
\e{cor}

We will later choose $\beta=4\theta_0+1$; see \eqref{eq:beta0defn}. In the remainder of this section, we prove Theorem \ref{thm:concentration}.



\subsection{Initial moment estimate}
In this section, we state an initial estimate which is the starting point of our analysis. It expresses the moments of $F_i(z)$ as sums over solutions of appropriate Diophantine conditions weighted by entries of the Green's function. We denote
\beq\label{eq:[N]defn}
[N]:=\{1,\ldots,N\}.
\eeq

Given an integer $p\geq 1$ and a vector $v\in \Z^d$, we define the set
\begin{equation} \label{eq:Lvdefn}
\begin{aligned}
    \mathcal{L}^{2p}_v:=
    \Bigg\{(&\ul{k},\ul{l})=(k_1,\ldots,k_{2p},l_1,\ldots,l_{2p})
    \in [N]^{2p}\times [N]^{2p}\, :\\
    & k_\al \ne l_\al,\, \forall 1\le \al \le 2p \,\textnormal{ and }\,
    \sum_{\al'=1}^{2p} (k_{\al'}^q - l_{\al'}^q)= v_q,\, \forall 1\le q \le d  \Bigg\}.
    \end{aligned}
\end{equation}

The following lemma is the only place where our Assumption \ref{ass:omega} on the random variables enters.

\begin{lm}\label{lm:omegaexp}
Let $\om\in \R/\Z$ be distributed according to $\rho(\om)\d\om$ with $\rho:\R/\Z\to [0,\infty)$ satisfying $\|\rho'\|_\infty\leq C$. Then there exists a constant $C'>0$ such that
$$
\l|\int_0^1 e[a\om]\rho(\om)\d\om\r|
\leq  \mathbbm 1_{a=0}+\frac{C'}{|a|} \mathbbm{1}_{a\neq 0},
$$
for any $a\in\Z$.
\end{lm}

\begin{proof}
The case $a=0$ is trivial and the case $a\neq 0$ follows from integration by parts and the assumption that $\|\rho'\|_\infty\leq C$.
\end{proof}

Given a vector $v\in \Z^d$, we define the function
\beq\label{eq:fdefn}
f(v):=\prod_{q=1}^d \l(\mathbbm 1_{v_q=0}+\frac{C'}{v_q} \mathbbm{1}_{|v_q|\neq 0}\r)
\eeq
with $C'$ given from Assumption \ref{ass:omega} and Lemma \ref{lm:omegaexp}.

The following estimate sets the stage for this section.

\begin{prop}[Initial moment estimate]
\label{prop:initial}
Let $p\geq 1$ be an integer. We have
\beq
|\mathbb E_i[F_i(z)^{2p}]|
\leq \frac{1}{N^{2p}}
\sum_{v\in \Z^d} f(v) \sum_{(\ul{k},\ul{l})\in \curly{L}^{2p}_v}
 |G^{(i)}_{k_1,l_1}|\ldots |G^{(i)}_{k_p,l_p}| 
\eeq
\end{prop}

\begin{proof}
We recall that for each $1\leq i\leq N$, we write $\mathbb E_i$ for the partial expectation with respect to the random variables $\{\om_{i,q}\}_{1\leq q\leq d}$. 
From \eqref{eq:Fidefn}, the fact that $G^{(i)}$ is independent of the $\{\om_{i,q}\}_{1\leq q\leq d}$ and Lemma \ref{lm:omegaexp} we find that
$$
\begin{aligned}
&\l|\mathbb E_i[F_i(z)^{2p}]\r|\\
&=
\frac{1}{N^{2p}}
\l|\sum_{\substack{1\leq k_1,l_1,\ldots,k_p,l_p\leq N:\\ k_\al \ne l_\al}}
G^{(i)}_{k_1,l_1}\ldots G^{(i)}_{k_p,l_p} \mathbb E_i \l[e\l[ \sum_{q=1}^d(k_1^q+\ldots k_p^q -l_1^q-\ldots -l_p^1)\omega_{i,q}\r] \r]
\r|\\
&\leq  \frac{1}{N^{2p}}
\sum_{v\in \Z^d} f(v) \sum_{(\ul{k},\ul{l})\in \curly{L}^{2p}_v}
 |G^{(i)}_{k_1,l_1}|\ldots |G^{(i)}_{k_p,l_p}| \\
\end{aligned}
$$
as claimed.
\end{proof}

\subsection{Cardinality bound for $\mathcal L_v^{2p}$}
We can bound the cardinality $|\mathcal L_v^{2p}|$ based on the 2015 breakthrough of Bourgain-Demeter-Guth \cite{BDG} that proved Vinogradov's Main Conjecture. Later on, we will decompose $\mathcal L_v^p$ into a good and bad set and this result this will play the role of an a priori estimate on the number of terms that need to be treated. The bound uses that $4p<d(d+1)$ which can be readily verified for our choices of $p=\floor{\frac{17d}{16}}$ and $d\geq 18$.

\begin{thm}[Cardinality bound \cite{BDG}]
\label{thm:vino}
For every $\eps>0$, there exists $C_\eps>1$ so that for all $N\geq 1$, it holds that 
\beq\label{eq:bdg}
|\mathcal L_v^{2p}|\leq C_\eps N^{2p+\eps}.
\eeq
\end{thm}

As a point of reference, we note that $N^{2p}$ is a trivial lower bound on $|\mathcal L_0^p|$ which can be seen by considering the diagonal solutions $k_q=j_q$ for all $1\leq q\leq d$, so \eqref{eq:bdg} is essentially sharp.

\be{rmk}
Here and in the following, we often suppress the dependence of various multiplicative constants on parameters such as $d$ and $p$ whenever this dependence plays no role for the ensuing argument. Moreover, the value of constants such as $C$ or $C_p$ may change from line to line.
\e{rmk}

\begin{proof}
Define the set
$$
\tilde{\mathcal{L}}^{2p}_v:=
    \Bigg\{(\ul{k},\ul{l})=(k_1,\ldots,k_{2p},l_1,\ldots,l_{2p})
    \in [N]^{2p}\times [N]^{2p}\, :\\
    \sum_{\al'=1}^{2p} (k_{\al'}^q - j_{\al'}^q)= v_q,\, \forall 1\le q \le d  \Bigg\}
$$
which differs from $\mathcal L_v^{2p}$ in that it can have $k_\al=l_\al$. We recall that $e[x]=\exp(2\pi i x)$. We have
$$
|\mathcal L_v^{2p}|
\leq |\tilde{\mathcal L}_v^{2p}|=\int_{[0,1]^d}
\l|\sum_{n=1}^N 
e\l[\sum_{q=1}^d \xi_q n^q\r]\r|^{4p} 
e\l[- \sum_{q=1}^d \xi_q v_q\r] \d \xi_1\ldots \d\xi_d,
$$
where the equality follows by expanding the power and using orthonormality of the Fourier basis. The triangle inequality then implies $|\tilde{\mathcal L}_v^{2p}|\leq |\tilde{\mathcal L}_0^{2p}|$. The key input is Theorem 1.1 in \cite{BDG} which says 
$$
|\tilde{\mathcal L}_0^{2p}|\leq C_\eps N^\eps (N^{2p}+N^{4p-\frac{d(d+1)}{2}})
$$
and the second term is subleading for $4p<d(d+1)$.
\end{proof}

\subsection{Cardinality bounds for the set of bad indices}
Taking a closer look at the right-hand side of Proposition \ref{prop:initial}, the task is clear: We need to control the size of the Green's function entries $|G^{(i)}_{k,l}|$ over the index set $\curly{L}_v^{2p}$ while retaining control over the size of the index set. While the bound from Theorem \ref{thm:vino} is helpful in this regard, it is far from sufficient because it yields no control on the size of the Green's function entries. Moreover, the bound \eqref{eq:bdg} is too course in general because it ignores the off-diagonal condition $k_\al\neq l_\al$ in the definition of the set $\curly{L}_v^{2p}$ (cf.\ the proof of Theorem \ref{thm:vino}), but we need to use this condition because it reflects the centering of the fluctuation variable $F_i(z)$, cf.\ \eqref{eq:Fidefn} without which there is no concentration.

To address this problem, we develop a refined pigeonhole principle which is at the heart of our proof. The upshot of these considerations is the bound below on the cardinality of the ``bad sets'' defined as follows.

\be{defn}[The bad sets]\label{defn:bad}
Let $r\geq 0$ be an integer and let $\gam>0$. Let
\begin{equation}
    \mathcal{B}_v^{r,\gam}:=\{(k_1,\ldots,k_{2p},l_1,\ldots,l_{2p}) \in \mathcal{L}^{2p}_v: |\{1\leq \al\leq 2p: |G^{(i)}_{k_\al,l_\al}| \le N^{-\gamma}\}| = r  \}.
\end{equation}
\e{defn}

In words, on the bad set $\mathcal{B}_v^{r,\gam}$ there are ``few'' ($=r$) Green's function entries that are ``small'' (less than $N^{-\gam}$).
The optimal choices for the parameters $r$ and $\gam$ will be determined later. 

We also define the counting function
\begin{equation}\label{eq:grNdefn}
    g_r(N):=N^{2r} N^{(2p-d_0 -r)(1+2\gam+2\theta)}, \qquad \textnormal{ with } d_0:=\floor{\frac{d}{2}}.
\end{equation}

The following cardinality bound for the bad sets is essential. 

\begin{thm}[Cardinality bound for the bad set $\curly{B}_v^{r,\gam}$]
\label{thm:bad}
There exists a constant $C_{p,r}>0$ so that for all $N\geq 1$ and all vectors $v\in \Z^d$,
\beq\label{eq:cardbad}
|\curly{B}_v^{r,\gam}|
\leq
C_{p,r}\l(N^{2d_0-2} N^r N^{(2p-r)(2\gam+2\theta)}+\max\{g_0(N),g_r(N)\}\r)
\eeq
\end{thm}

Theorem \ref{thm:bad}  quantifies the extent to which the cardinality of the bad set is subleading compared to the size of the whole index set $\mathcal{L}^{2p}_v$, i.e., compared to $N^{2p+\eps}$ according to Theorem \ref{thm:vino} and so it makes precise the notion that ``bad indices are rare''. We defer the proof of this important estimate to the next section. 

\subsection{A priori estimate on the Green's function}
The Ward identity implies weak a priori estimates on the size of the $|G^{(i)}_{k,l}|$ defined in \eqref{eq:Gidefn}.

\begin{lm} \label{lm:BoundonG}
Let $z\in \C$ with $\Im[z]= N^{-\theta}$. Then
\begin{equation}
\label{eq:BoundonG}
|G_{k,l}^{(i)}| \le N^{\theta},\qquad  \forall 1\leq k,l\leq N.
\end{equation}
\end{lm}
\begin{proof}
Denote $G^{(i)}=G$. Note that $G=\frac{1}{A-z}$ with $A$ a Hermitian $N\times N$ matrix. From the Ward identity and $\Im[G_{k,k}]\leq|G_{k,k}|\leq \|G\|\leq |z|^{-1}\leq N^{\theta}$, we have
\begin{equation}
    \sum_{l=1}^{N} |G_{k,l}|^2 = \frac{\Im[G_{k,k}]}{\Im[z]} \le N^{2\theta}
\end{equation}
and this implies \eqref{eq:BoundonG}.
\end{proof}

\subsection{Choice of parameters and conclusion}
Assuming Theorem \ref{thm:bad} holds, we can now prove Theorem \ref{thm:concentration} by choosing near-optimal parameters $p,r,\gam$ (depending on the model parameter $d$). 

\begin{proof}[Proof of Theorem \ref{thm:concentration}]
We apply Proposition \ref{prop:initial} noting that $\curly{L}_v^{2p}\neq\varnothing$ implies 
$$
|v_q|\leq C_{p,q}N^q\leq C_p N^d,\qquad \forall 1\leq q\leq d
$$
so that
$$
|\mathbb E_i[F_i(z)^{2p}]|
\leq \frac{1}{N^{2p}}
\sum_{\substack{v\in \Z^d:\\ |v_q|\leq C_pN^d}} f(v) \sum_{(\ul{k},\ul{l})\in \curly{L}^{2p}_v}
 |G^{(i)}_{k_1,l_1}|\ldots |G^{(i)}_{k_p,l_p}|
$$
Next, we decompose the index set $\curly{L}_v^{2p}$ as follows.
$$
\curly{L}_v^{2p}=\bigcup_{r'=0}^r \mathcal{B}_v^{r',\gam} 
\cup \curly{G}_v^r,
\qquad \textnormal{ with } \curly{G}_v^r:=\curly{L}_v^{2p}\setminus \bigcup_{r'=0}^r \mathcal{B}_v^{r',\gam}.
$$
We call $\curly{G}_v^r$ the ``good set''. We implement this decomposition to estimate the Green's function. By Definition \ref{defn:bad}, the good set contains at least $r+1$ ``small'' Green's function elements $|G^{(i)}_{k,l}|\leq N^{-\gam}$. The remaining Green's function entries on the good set are bounded by the a priori estimate $N^\theta$ from Lemma \ref{lm:BoundonG}. Applying analogous bounds for the various bad sets, we obtain
\beq\label{eq:goodbad}
\begin{aligned}
&|\mathbb E_i[F_i(z)^{2p}]|\leq \frac{1}{N^{2p}}
\sum_{\substack{v\in \Z^d:\\ |v_q|\leq C_pN^d}} f(v) \sum_{(\ul{k},\ul{l})\in \curly{L}^{2p}_v}
 |G^{(i)}_{k_1,l_1}|\ldots |G^{(i)}_{k_p,l_p}|\\ 
&= \frac{1}{N^{2p}}
\sum_{\substack{v\in \Z^d:\\ |v_q|\leq C_pN^d}} f(v) 
\l(
\sum_{(\ul{k},\ul{l})\in \curly{G}_v^r}
  |G^{(i)}_{k_1,l_1}|\ldots |G^{(i)}_{k_p,l_p}| +\sum_{r'=0}^r
\sum_{(\ul{k},\ul{l})\in \curly{B}_v^{r',\gam}}
 |G^{(i)}_{k_1,l_1}|\ldots |G^{(i)}_{k_p,l_p}| 
\r)\\
  &\leq 
  \frac{1}{N^{2p}}
\sum_{\substack{v\in \Z^d:\\ |v_q|\leq C_pN^d}} f(v) 
\l(
  N^{-(r+1)\gam}N^{ (2 p-r-1)\theta}|\curly{G}_v^r| +\sum_{r'=0}^r N^{-r'\gam} N^{(2p-r')\theta} |\curly{B}_v^{r',\gam}|
 \r)\\
  &\leq 
  \frac{C_{p,r} C_\eps}{N^{2p}}
\l(\sum_{\substack{v\in \Z^d:\\ |v_q|\leq C_pN^d}} f(v) \r)
\l( N^{-(r+1)\gam}N^{ (2 p-r-1)\theta}N^{ 2p+\eps}+\sum_{r'=0}^r N^{-r'\gam} N^{(2p-r')\theta}\Phi_{r'}(N)\r),
 \end{aligned}
 \eeq
 where the last step uses Theorem \ref{thm:vino} (together with the trivial estimate $|\curly{G}_v^r|\leq |\curly{L}^{2p}_v|$) and Theorem \ref{thm:bad}.  Here we introduced the function
$$
\Phi_{r'}(N)=N^{2d_0-2} N^{r'} N^{(2p-r')(2\gam+2\theta)}+\max\{g_0(N),g_{r'}(N)\}
$$ 

Before we analyze the exponents of $N$ further, we note that the sum over $v$ can now be performed. Recalling Definition \eqref{eq:fdefn} of $f$, we have
\beq\label{eq:festimate}
\sum_{\substack{v\in \Z^d:\\ |v_q|\leq C_pN^d}} f(v)
=\prod_{q=1}^d \l(1+C'\sum_{v_q=1}^{C_pN^p} \frac{1}{|v_q|}\r)
\leq C_p (\log N)^d.
\eeq
We see that this term is of logarithmic size and thus (almost) irrelevant.

We recall that the claim of Theorem \ref{thm:concentration} is to estimate $|\mathbb E_i[F_i(z)^{2p}]|$ by $(N^{-\frac{1}{36}+\eps'}/\Im[z])^{2p}=N^{2p(\theta-\frac{1}{36}+\eps')}$ for any $\eps'>0$. From \eqref{eq:goodbad} and \eqref{eq:festimate}, we see that the task is to show that
\beq\label{eq:mainremains}
C_{\eps,p,r} (\log N)^d N^{-2p(\theta+1)}\l(N^{-(r+1)\gam+ (2 p-r-1)\theta+2p+\eps}+\sum_{r'=0}^r N^{-r'\gam} N^{(2p-r')\theta}\Phi_{r'}(N)\r)
\leq N^{-\frac{2p}{36}+2p\eps'}
\eeq
for sufficiently large $N$. We can still choose the parameters $\theta,r,\gam$. Simplifying the left-hand side in \eqref{eq:mainremains} gives
$$
\begin{aligned}
& (\log N)^dN^{-2p(\theta+1)}\l(N^{-(r+1)\gam+ (2 p-r-1)\theta+2p+\eps}+\sum_{r'=0}^r N^{-r'\gam} N^{(2p-r')\theta}\Phi_{r'}(N)\r)\\
&=(\log N)^d \Big(N^{-(r+1)\gam+ -(r+1)\theta+\eps}+
N^{2d_0-2-2p+4p(\gam+\theta)}
\sum_{r'=0}^r N^{r'(1-3\gam-3\theta)}\\
&\qquad+N^{(2p-d_0)(1+2\gam+2\theta)}
\sum_{r'=0}^r N^{-r'(\gamma+\theta)} \max\{1,N^{r'(1-2\gam-2\theta)}\}\Big)\\
&\leq C(\log N)^d(N^{x_1}+N^{x_2}+N^{x_3}). 
\end{aligned}
$$
The last step uses that the summands are either monotonically increasing or decreasing in $r'$ (depending on the value of $\gam+\theta$) and introduces the three exponents
\beq\label{eq:xi}
\begin{aligned}
x_1=&-(r+1)(\gam+\theta)+\eps,\\
x_2=&-2p+2d_0-2+4p(\gam+\theta)+r(1-3\gam-3\theta)_+\\
x_3=&-2p+(2p-d_0)(1+2\gam+2\theta)+r((1-2\gam-2\theta)_+-\gam-\theta)_+,
\end{aligned}
\eeq
with $(y)_+=\max\{y,0\}$ denoting the positive part of a real number $y$.

In view of the various case distinctions we impose 
$$
\gam+\theta=:\theta'\leq \frac{1}{3},
$$
 in which case the exponents simplify to
\beq\label{eq:xisimplified}
\begin{aligned}
x_1=&-(r+1)\theta'+\eps,\\
x_2=&-2p+2d_0-2+4p\theta'+r(1-3\theta')\\
x_3=&2p\theta'-d_0(1+2\theta')+r(1-3\theta'),
\end{aligned}
\eeq
Elementary estimates show that the parameter values
\beq\label{eq:choice}
\theta'=\frac{1}{5},\qquad r=\frac{22d}{51}, \qquad p=\floor{\frac{17d}{16}}
\eeq
gives $\max\{x_1,x_2,x_3\}<-\frac{2p}{36}-2p\eps'$ for all $d\geq 18$ and all $\eps'>0$ provided that $\eps$ is chosen sufficiently small. (We arrived at the choice  \eqref{eq:choice} by assuming that $p$ and $r$ are linear multiples of $d$ and optimizing $x_1,x_2,x_3$ among this class.) This proves Theorem \ref{thm:concentration}.
\end{proof}

\section{Proof of Theorem \ref{thm:bad}}
At the heart of our proof is a somewhat delicate pigeonhole argument which rests on structural aspects of the Ward identity and a conditional uniqueness result for the Diophantine equations defining $\curly{L}_v^{2p}$ that follows from the Newton-Girard identities.

\subsection{The Ward identity within rows and columns}
A fundamental observation that guides our approach is that the a priori bound in Lemma \ref{lm:BoundonG} (a direct consequence of the Ward identity) can be improved if many of the Green's function entries lie in the same row (or column). This is made precise by the following lemma.

\begin{lm}\label{lm:ap}
Let $z\in\C$ with $\Im[z]= N^{-\theta}$ and let $1\leq i\leq N$. We have 
\beq
\label{eq:apfull}
\l|\setof{1\leq k,l\leq N}{|G^{(i)}_{k,l}|> N^{-\gam}}\r|\leq N^{1+2\gam+2\theta}
\eeq
and for every $k_0,l_0\in\{1,\ldots,N\}$,
\beq
\label{eq:aprow}
\begin{aligned}
\l|\setof{1\leq l\leq N}{|G^{(i)}_{k_0,l}|> N^{-\gam}}\r|\leq N^{2\gam+2\theta},\\
\l|\setof{1\leq k\leq N}{|G^{(i)}_{k,l_0}|> N^{-\gam}}\r|\leq N^{2\gam+2\theta}.
\end{aligned}
\eeq
\end{lm}

\begin{proof}
    We denote $G=G^{(i)}$. Due to the Ward identity and symmetry, we have for each $k_0$,
    \begin{equation}\label{eq:wardrow}
        \sum_{l=1}^N |G_{k_0,l}|^2 =  \frac{\Im[G_{k_0,k_0}]}{\Im[z]} \le N^{2\theta}.
    \end{equation}
    This implies the first bound in \eqref{eq:aprow}, while the second one follows from symmetry of $G$. Finally, \eqref{eq:apfull} follows from summing \eqref{eq:wardrow} over $1\leq k_0\leq N$ and so Lemma \ref{lm:ap} is proved.
\end{proof}
The fundamental question we thus need to investigate next is how the structural property of lying in the same row interacts with the Diophantine conditions that define the index set $\curly{L}_v^{2p}$ in \ \eqref{eq:Lvdefn}. 

\subsection{A dichotomy for the index pairs}
We note a simple dichotomy: Either many index pairs lie in the same row (or column) or many do not. The precise version is given in the following lemma, which involves an integer parameter $1\leq s\leq 2p$ (which is later chosen as $s=d_0=\floor{d/2}$, so relatively large).

\begin{lm}[Index set dichotomy] \label{lm:Disj}
Consider the collection of pairs
$$
\curly{P}=\setof{(k_\al,l_\al)\in \Z\times\Z}{1\leq \al\leq 2p,\, k_\al\neq l_\al}.
$$
For every integer $s\geq 1$, one of the following two statements holds.
\begin{itemize}
\item[(a)]  There are $\nu_1,\ldots,\nu_{2s-2}\in\Z$ such that
for all $1\leq \al\leq 2p$, either $k_\al$ or $l_\al$ lie in $\{\nu_1,\ldots,\nu_{2s-2}\}$.
\item[(b)] There exist distinct $1\leq \al_1,\ldots,\al_s\leq 2p$ such that
$$
\{k_{\al_1},\ldots,k_{\al_s}\}\cap \{l_{\al_1},\ldots,l_{\al_s}\}=\varnothing,
$$
\e{itemize}
\end{lm}

\begin{proof}
We induct in $s$. The base case $s=1$ is trivial since (b) holds by assumption. 

For the induction step, assume that the claim holds for $s-1$.
If case (a) occurred for $s-1$, then it also occurs for $s$ (it is a monotone condition), so we may assume that case (b) occurs for $s-1$. That is, there exist 
distinct $1\leq \al_1,\ldots,\al_{s-1}\leq 2p$ such that
\beq\label{eq:emptycond}
\{k_{\al_1},\ldots,k_{\al_{s-1}}\}\cap \{l_{\al_1},\ldots,l_{\al_{s-1}}\}=\varnothing.
\eeq
Consider the remaining pairs $(k_{\al'},l_{\al'}) \in\curly{P}$ with $\al'\neq \al_1,\ldots,\al_{s-1}$. We distinguish two cases. The first case is that we can find $\al'\neq \al_1,\ldots,\al_{s-1}$ such that (b) holds with $\al_s=\al'$, in which case the induction step is completed. Otherwise, we have
\beq\label{eq:nonemptycond}
\{k_{\al_1},\ldots,k_{\al_{s-1}},k_{\al'}\}\cap \{l_{\al_1},\ldots,l_{\al_{s-1}},l_{\al'}\}\neq \varnothing,\qquad  \forall \al'\neq \al_1,\ldots,\al_{s-1}.
\eeq
We set
$$
\{\nu_1,\ldots,\nu_{2s-2}\}=
\{k_{\al_1},\ldots,k_{\al_{s-1}},l_{\al_1},\ldots,l_{\al_{s-1}}\}
$$
and note that condition (a) now follows from \eqref{eq:emptycond}, \eqref{eq:nonemptycond} and $k_{\al'}\neq l_{\al'}$. This proves Lemma \ref{lm:Disj}.
\end{proof}

\subsection{Algebraic conditional uniqueness for Diophantine equations}
According to Lemma \ref{lm:Disj}, there are either $2s-2$ indices in the same row or column as in case (a) or $s$ index pairs are from different rows/columns as in case (b). In case (a), the refined a priori bound \eqref{eq:aprow} from Lemma \ref{lm:ap} is helpful, so case (b) has to be understood next. 

Here, we now show that we get case (b) is suppressed for algebraic reasons, namely by the Diophantine equations that constrain the index set $\curly{L}_v^p$. This is made precise in Lemma \ref{lm:uniq} below, a purely algebraic conditional uniqueness result which we observe here but which we suspect is well-known to experts in number theory.

\begin{lm}[Conditional uniqueness condition for Diophantine equations] \label{lm:uniq}
Let $n$ be an integer and consider two collections $\curly{P}_1$ and $\curly{P}_2$ of integer pairs
$$
\curly{P}_1=\setof{(i_\al,j_\al)\in \Z\times \Z}{1\leq \al\leq n},\qquad 
\curly{P}_2=\setof{(i_\al',j_\al')\in \Z\times \Z}{1\leq \al\leq n}
$$
subject to the disjointness conditions
\beq\label{eq:disj}
\{i_1,\ldots,i_n\}\cap \{j_1,\ldots,j_n\} =\varnothing,\qquad \{i_1',\ldots,i_n'\}\cap \{j_1',\ldots,j_n'\} =\varnothing.
\eeq
Suppose that there exist $d_1,\ldots,d_{2n}\in \Z$ such that
\beq\label{eq:ParVino}
\sum_{\al=1}^n (j_\al^q-i_\al^q)=\sum_{\al=1}^n ((j_\al')^q-(i_\al')^q)=d_q,\qquad \forall 1\leq q\leq 2n.
\eeq
Then
\beq\label{eq:dioclaim}
\{i_1,\ldots,i_n\}=\{i_1',\ldots,i_n'\}
\textnormal{ and }
\{j_1,\ldots,j_n\}=\{j_1',\ldots,j_n'\}
\eeq
\e{lm}

\be{proof}
The set of Diophantine equations \eqref{eq:ParVino} imply the following analog without the negative sign,
$$
\sum_{\al=1}^n((i_\al')^q+j_\al^q)=\sum_{\al=1}^n(i_\al^q+(j_\al')^q),\qquad \forall 1\leq q\leq 2n.
$$
By the Newton-Girard identities, the power sums up to order $2n$ determine the corresponding elementary symmetric polynomials up to the same order. Consequently, we have the equality of the two polynomials 
$$
\prod_{\al=1}^n (x-i'_\al)(x-j_\al)=\prod_{\al=1}^n (x-i_\al)(x-j'_\al).
$$
The equality of the polynomials implies the equality of their root sets, i.e., 
$$
\{i'_1,\ldots i'_n,j_1,\ldots,j_n\}=\{i_1,\ldots i_n,j'_1,\ldots,j_n\}.
$$
The claim \eqref{eq:dioclaim} now follows from the disjointness assumption \eqref{eq:disj}.
\e{proof}

\subsection{Proof of Theorem \ref{thm:bad}}
We now have all the tools in hand to prove Theorem \ref{thm:bad} via a refined pigeonhole principle.

\begin{proof}[Proof of Theorem \ref{thm:bad}]
We fix integers $p,r,N\geq 1$, a number $\gam>0$ and an arbitrary vector $v\in \Z^d$. Given any index list $(\ul{k},\ul{l})=(k_1,\ldots,k_{2p},l_1,\ldots,l_{2p})
    \in \curly{L}_v^{2p}$, we apply Lemma \ref{lm:Disj} with 
    $
    s=d_0=\floor{d/2}
    $
     to the corresponding list of pairs $\curly{P}$. (This choice of $s$ turns out to be optimal for Step 1 in the proof of Lemma \ref{lm:caseb} later on.) Hence, we can decompose $\curly{L}_v^p$, and consequently the bad set $\mathcal{B}=\mathcal{B}_v^{p,r,\gam}$, as follows
\beq\label{eq:baddecomp}
\curly{B}=\curly{B}_{(a)}\cup \curly{B}_{(b)}
\eeq
where $\curly{B}_{(x)}$ with $x\in\{a,b\}$ is the set of $(\ul{k},\ul{l})\in \curly{B}$ such that the respective case occurs in Lemma \ref{lm:Disj}.

The following two lemmas bound the cardinalities of $\curly{B}_{(a)}$ and $\curly{B}_{(b)}$.

\be{lm}[Case (a) bound]\label{lm:casea}
There exists a constant $C_{p,r}>0$ such that for all $N\geq1$,
\beq
|\curly{B}_{(a)}|\leq C_{p,r} N^{2d_0-2} N^{r }N^{(2p-r) (2\gam+2\theta)}
\eeq
\e{lm}

\be{lm}[Case (b)]\label{lm:caseb}
There exists a constant $C_{p,r}>0$ such that for all $N\geq1$,
\beq
|\curly{B}_{(b)}|\leq C_{p,r} \max\{g_0(N),g_r(N)\}.
\eeq
\e{lm}

Considering \eqref{eq:baddecomp}, we see that Theorem \ref{thm:bad} follows from Lemmas \ref{lm:casea} and \ref{lm:caseb}.
\e{proof}

\subsection{Proof of Lemmas \ref{lm:casea} and \ref{lm:caseb}}
\be{proof}[Proof of Lemma \ref{lm:casea}]
We consider the constraints that exist on a generic element $(\ul{k},\ul{l})\in\curly{B}_{(a)}$ and use this to estimate $|\curly{B}_{(a)}|$ through basic combinatorics. Since we are in case (a) of Lemma \ref{lm:Disj} with $s=d_0$, there exist $\nu_1,\ldots,\nu_{2d_0-2}\in [N] =\{1,\ldots,N\}$, such that for all $1\leq \al\leq 2p$, either $k_\al$ or $l_\al$ lie in $\{\nu_1,\ldots,\nu_{2d_0-2}\}$. It will be convenient to introduce
$$
m_\al=
\be{cases}
k_\al,\qquad \textnormal{if } k_\al\in \{\nu_1,\ldots,\nu_{2d_0-2}\},\\
l_\al,\qquad \textnormal{otherwise.}
\e{cases}
$$
and
$$
n_\al=
\be{cases}
l_\al,\qquad \textnormal{if } k_\al\in \{\nu_1,\ldots,\nu_{2d_0-2}\},\\
k_\al,\qquad \textnormal{otherwise.}
\e{cases}
$$
First, we note that there are at most $N^{2d_0-2}$ ways to choose the $\nu_1,\ldots,\nu_{2d_0-2}\in [N]$. Second, we choose for every $1\leq \al\leq 2p$ whether $k_\al\in \{\nu_1\ldots,\nu_{2d_0-2}\}$ or not and for this there are $2^{2p}$ options. After this step, it is determined whether $(k_\al,l_\al)=(m_\al,n_\al)$ or $(n_\al,m_\al)$ and so it remains to count the options for $m_\al$ and $n_\al$.

Regarding the number of choices for $m_1,\ldots,m_{2p}$, we note since case (a) of Lemma \ref{lm:Disj} applies, we have $m_1,\ldots,m_{2p}\in \{\nu_1,\ldots,\nu_{2d_0-2}\}$ and so there are at most $(2d_0-2)^{2p}$ choices. 

To summarize the considerations so far, we have the combinatorial factor 
\beq\label{eq:sofar}
2^{2p}(2d_0-2)^{2p}N^{2d_0-2}
\eeq
which accounts for the number of choices of everything except the $n_1,\ldots,n_{2p}$. For these, we shall use that configurations $(\ul{k},\ul{l})\in\curly{B}_{(a)}\subset \curly{B}$ are constrained further because they must belong to the bad set. Indeed, recalling the Definition \ref{defn:bad} of the bad set, there must be $2p-r$ choices of $\al$ so that $|G^{(i)}_{k_\al,l_\al}|>N^{-\gam}$ is ``large''. A simple but important observation is that
$$
\big\{|G^{(i)}_{k_\al,l_\al}|\,:\, 1\leq \al\leq 2p\big\}=\setof{|G^{(i)}_{m_\al,n_\al}|}{1\leq \al\leq 2p}
$$
because $G^{(i)}$ is a symmetric matrix. Consequently, the collection $\{|G^{(i)}_{m_\al,n_\al}|\}_{1\leq \al\leq 2p}$ must also contain $2p-r$ large elements. The advantage of this collection is that the $|G^{(i)}_{m_\al,n_\al}|$ all belong to the same $2d_0-2$ rows since $m_1,\ldots,m_{2p}\in \{\nu_1,\ldots,\nu_{2s-2}\}$ and every row contains at most $N^{2\gam+2\theta}$ large entries by the second part of Lemma \ref{lm:ap}.

 Hence, the number of choices for $n_1,\ldots,n_{2p}$ can be estimated as follows: First we choose which $2p-r$ of the $1\leq \al\leq 2p$ correspond to large Green's function entries and for this there are $\binom{2p}{r}$ options. Second, we distribute the $2p-r$ of the $n_\al$'s corresponding to large $|G^{(i)}_{m_\al,n_\al}|$ among the at most $N^{2\gam+2\theta}$ options in their assigned row $m_\al$, resulting in a total number of $N^{(2p-r)(2\gam+2\theta)}$ options. Third, we distribute the remaining $r$ of the $n_\al$'s among the at most $N$ options in their assigned row $m_\al$ and for this there are at most $N^r$ options. In summary, we have shown that the number of choices for the $n_1,\ldots,n_{2p}$ is bounded by
$$
\binom{2p}{r}N^{(2p-r)(2\gam+2\theta)} N^r
$$
(As a point of reference, a naive estimate on the choices for $n_1,\ldots,n_{2p}$ is of course $N^{2p}$.) 

Multiplying this by the other combinatorial factor from \eqref{eq:sofar} gives the bound
$$
|\curly{B}_{(a)}|\leq \binom{2p}{r} 2^{2p}(2d_0-2)^{2p}N^{2d_0-2}N^{(2p-r)(2\gam+2\theta)}N^r.
$$
This proves Lemma \ref{lm:casea} with the constant $C_{p,r}=\binom{2p}{r} 2^{2p}(2d_0-2)^{2p}$.
\e{proof}

\be{proof}[Proof of Lemma \ref{lm:caseb}]
This proof is one of the essential technical parts of our argument and it uses the Lemmas \ref{lm:ap}, \ref{lm:Disj}, and \ref{lm:uniq} that were established earlier. It will be convenient to denote
$$
\ul{k}=(k_1,\ldots,k_{2p})=(\mathbf{k},k_{2p-d_0},\ldots,k_{2p})\quad \textnormal{ with } \quad \mathbf{k}=(k_1,\ldots,k_{2p-d_0})
$$
Given $\mathbf{k},\mathbf{l}\in [N]^{2p-d_0}$, we define the following subsets of $\curly{B}_{(b)}$ labeled by the first $2p-d_0$ elements of each sequence.
\begin{equation}
\begin{aligned}
    \mathcal{B}_{(b)}(\mathbf{k},\mathbf{l})
    :=\big\{(&\mathbf{k},k_{2p-d_0+1},\ldots,k_{2p},\mathbf{l},l_{2p-d_0+1},\ldots,l_{2p}) \in \mathcal{B}_{(b)}\, :\,\\
      &\{k_{2p-d_0+1},\ldots,k_{2p}\} \cap \{ l_{2p-s+1},\ldots,l_{2p}\} = \varnothing \}  \big\}
\end{aligned}
\end{equation}
We can use these sets to further subdivide the set $\curly{B}_{(b)}$ because we know that every element $(k_1,\ldots,k_{2p},l_1,\ldots,l_{2p}) \in \mathcal{B}_{(b)}$ satisfies case (b) in Lemma \ref{lm:Disj}. That is, modulo permutation, the last $d_0$ elements of the sequences are disjoint and so a union bound gives
\beq\label{eq:permutations}
|\curly{B}_{(b)}|
\leq (2p)! \l|\bigcup_{\mathbf{k},\mathbf{l}\in[N]^{2p-d_0}}\mathcal{B}_{(b)}(\mathbf{k},\mathbf{l})\r|
\eeq
Below, we prove that
\beq\label{eq:cardclaim}
\l|\bigcup_{\mathbf{k},\mathbf{l}\in[N]^{2p-d_0}}\mathcal{B}_{(b)}(\mathbf{k},\mathbf{l}),
\r|
\leq C_{p,r} \max\{g_0(N),g_r(N)\}.
\eeq
which together with \eqref{eq:permutations} implies Lemma \ref{lm:caseb}.

Thus, it remains to prove \eqref{eq:cardclaim}. This is done in two steps.\\

\dashuline{Step 1.} We first estimate the size of each individual set appearing in \eqref{eq:cardclaim},
\beq
\label{eq:claim1}
|\mathcal{B}_{(b)}(\mathbf{k},\mathbf{l})|\leq (d_0!)^2,
\qquad \forall \mathbf{k},\mathbf{l}\in [N]^{2p-d_0}.
\eeq
We will derive this from the following exact characterization of these sets. It says that if $\mathcal{B}_{(b)}(\mathbf{k},\mathbf{l})\neq \varnothing$, then it is generated by permuting the last entries. More formally, if there exists some $(\mathbf{k},k_{2p-d_0},\ldots,k_{2p},\mathbf{l},l_{2p-d_0},\ldots,l_{2p}) \in \curly{B}_{(b)}(\mathbf{k},\mathbf{l})$, then
$\mathcal{B}_{(b)}(\mathbf{k},\mathbf{l})$ is in fact equal to the set of all $(\mathbf{k},k'_{2p-d_0},\ldots,k'_{2p},\mathbf{l},l'_{2p-d_0},\ldots,l'_{2p})$
satisfying
\beq\label{eq:aim}
\{k_{2p-d_0+1},\ldots,k_{2p}\}=\{k_{2p-d_0+1}',\ldots,k_{2p}'\},\qquad \{l_{2p-d_0+1},\ldots,l_{2p}\}=\{l_{2p-d_0+1}',\ldots,l_{2p}'\}.
\eeq
Note that this characterization implies \eqref{eq:claim1} since the number of non-trivial permutations of the last $d_0$ elements is bounded by $d_0!$. 

To complete Step 1, we need to prove this exact characterization of $\mathcal{B}_{(b)}(\mathbf{k},\mathbf{l})$. This part uses Lemma \ref{lm:uniq}. Consider two elements
$$
(\mathbf{k},k_{2p-d_0},\ldots,k_{2p},\mathbf{l},l_{2p-d_0},\ldots,l_{2p}),(\mathbf{k},k'_{2p-d_0},\ldots,k'_{2p},\mathbf{l},l'_{2p-d_0},\ldots,l'_{2p}) \in \curly{B}_{(b)}(\mathbf{k},\mathbf{l})
$$
for which we aim to prove \eqref{eq:aim}. By definition of $\mathcal{B}_{(b)}(\mathbf{k},\mathbf{l})$, we have
\beq\label{eq:verif1}
\{k_{2p-d_0+1},\ldots,k_{2p}\} \cap \{ l_{2p-d_0+1},\ldots,l_{2p}\} =\{k'_{2p-d_0+1},\ldots,k'_{2p}\} \cap \{ l'_{2p-d_0+1},\ldots,l'_{2p}\} = \varnothing.
\eeq
Now we recall that elements of $\curly{B}_{(b)}$ also lie in $\curly{L}_v^p$ and therefore solve the Diophantine equations from Definition \eqref{eq:Lvdefn} of $\curly{L}_v^p$. Solving the resulting sets of equations for the last $d_0$ elements, we obtain
\beq\label{eq:verif2}
\begin{aligned}
&\sum_{\al=1}^{d_0} (k_{2p-d_0+\al}^q-l_{2p-d_0+\al}^q)=\sum_{\al=1}^{d_0} ((k'_{2p-d_0+\al})^q-(l'_{2p-d_0+\al})^q)\\
&=d_q:=v_q-\sum_{{\al'}=1}^{2p-d_0}(k_{\al'}^q - l_{\al'}^q),\qquad \forall{1\leq q\leq d}.
\end{aligned}
\eeq
We see that \eqref{eq:verif1} and \eqref{eq:verif2} verify the conditions of Lemma \ref{lm:uniq} if we set $(i_\al,j_\al)=(k_{2p-d_0+\al},l_{2p-d_0+\al})$ and $(i'_\al,j'_\al)=(k'_{2p-d_0+\al},l'_{2p-d_0+\al})$ and choose $n=d_0=\floor{d/2}$ (noting also that $2n=2d_0\leq d$ as required). The conclusion of Lemma \ref{lm:uniq} is precisely \eqref{eq:aim}. This finishes Step 1.\\

\dashuline{Step 2.} While Step 1 controls the cardinality of an individual set $\curly{B}_{(b)}(\mathbf{k},\mathbf{l})$, we also need a bound on the number of terms in the union appearing in \eqref{eq:cardclaim}. (The trivial bound $N^{2(2p-d_0)}$ is insufficient for our purposes.)

In Step 2, we control the number of non-trivial choices of $\mathbf{k},\mathbf{l}\in [N]^{2p-d_0}$ that can lead to a non-empty $\curly{B}_{(b)}(\mathbf{k},\mathbf{l})$, i.e., we show
\beq
\label{eq:claim2}
\l|\setof{\mathbf{k},\mathbf{l}\in [N]^{2p-d_0}}{\curly{B}_{(b)}(\mathbf{k},\mathbf{l})\neq \varnothing}\r|
\leq C_{p,r}\max\{g_0(N),g_r(N)\}
\eeq  

This argument is of similar combinatorial flavor as the proof of Lemma \ref{lm:casea}, i.e., we estimate the cardinality in \eqref{eq:claim2} by studying the constraints on a generic pair $\mathbf{k},\mathbf{l}\in [N]^{2p-d_0}$ with $\curly{B}_{(b)}(\mathbf{k},\mathbf{l})\neq \varnothing$. We first recall that according to Definition \ref{defn:bad} elements of the bad set (of which $\curly{B}_{(b)}$ is a subset) hold $r$ ``small'' Green's function entries $|G^{(i)}_{k,l}|\leq N^{-\gam}$. In order to have $\curly{B}_{(b)}(\mathbf{k},\mathbf{l})\neq \varnothing$, there can be at most $r$ among the $(k_1,l_1),\ldots,(k_{2p-d_0},l_{2p-d_0})$ whose Green's function entries satisfy $|G^{(i)}_{k_\al,l_\al}|\leq N^{-\gam}$. Write $0\leq \varrho\leq r$ for the number of small entries. Given a value of $\varrho$, by permutation invariance, there are $\binom{2p-d_0}{\varrho}$ choices of the $\varrho$ indices among the $0\leq \al\leq 2p-d_0$ whose Green's function entries are small. 

Once we have selected the $\varrho$ indices among the $1\leq \al\leq 2p$ which have small Green's function entries, we have also fixed the remaining $2p-d_0-\varrho$ indices among the $1\leq \al\leq 2p-d_0$ which have large Green's function entries. For the first kind, there are trivially at most $N^2$ options for each $(k_\al,l_\al)$, resulting in a total of at most $N^{2\varrho}$ options. For the second kind, we note that the first part of Lemma \ref{lm:ap} implies there exist a total of at most $N^{1+2\gam+2\theta}$ index pairs $(k,l)\in [N]\times [N]$ for which $|G^{(i)}_{k,l}|$ is large. Hence, there are at most $N^{1+2\gam+2\theta}$ options for each $(k_\al,l_\al)$, resulting in a total of at most $N^{(2p-d_0-\varrho)(1+2\gam+2\theta)}$ options.

Altogether, taking into account the case distinction for the value of $0\leq \varrho\leq r$, these combinatorial considerations imply the estimate
$$
\begin{aligned}
\l|\setof{\mathbf{k},\mathbf{l}\in [N]^{2p-d_0}}{\curly{B}_{(b)}(\mathbf{k},\mathbf{l})\neq \varnothing}\r|
&\leq \sum_{\varrho=0}^r \binom{2p-d_0}{\varrho}N^{2\varrho}N^{(2p-d_0-\varrho)(1+2\gam+2\theta)}\\
&\leq \sum_{\varrho=0}^r \binom{2p-d_0}{\varrho}g_\varrho(N)\\
&\leq C_{p,r}\max\{g_0(N),g_r(N)\},
\end{aligned}
$$
where the last estimate uses that $r\mapsto g_r(N)$ is either monotonically increasing or monotonically decreasing. This proves \eqref{eq:claim2} and thus completes Step 2.

Finally, we note that the estimates \eqref{eq:claim1} and \eqref{eq:claim2} proved in Steps 1 and 2 together imply \eqref{eq:cardclaim} via the union bound. This completes the proof of Lemma \ref{lm:caseb}.
\e{proof}

\section{Step 4: Initial-scale estimate for large $\eta$}

\subsection{Stability analysis}
Our goal in Theorem \ref{thm:main1} is to show that $m_N(z)$ is close to $m_{\mathrm{\mathrm{MP}}}(z)$, the unique solution with positive imaginary part to \eqref{eq:mpquadratic}. This quadratic equation can be rearranged to
\beq\label{eq:mpquadratic'}
m_{\mathrm{\mathrm{MP}}}(z)+\frac{1}{z+zm_{\mathrm{\mathrm{MP}}}(z)}=0
\eeq

Equation \eqref{eq:step2} (and the control on the error term $|\curly{E}_i(z)|$ through Lemma \ref{lm:interlace} and Corollary \ref{cor:ld}) indicate that $m_{N}(z)$ satisfies an approximate version of this equation. It is essential for the proof that equation \eqref{eq:mpquadratic'} is stable in the sense that approximate solutions (with positive imaginary part) are close to $m_{\mathrm{\mathrm{MP}}}(z)$. 

\be{lm}[Stability]\label{lm:stability}
Let $z\in\curly{D}$. Suppose that $m$ satisfies
\beq\label{eq:stabilityass}
\l|m+\frac{1}{z+zm}\r|\leq \de
\eeq
for some $\de\leq 1$. Then
$$
\min\left\{|m - m_{\mathrm{\mathrm{MP}}}(z)|,  \l|m - \frac{1}{zm_{\mathrm{\mathrm{MP}}}(z)}\r|\right\} \le C \frac{\delta}{\kappa}.
$$
\e{lm}

Here we show that this stability follows by a simple substitution from the more widely known stability of the quadratic equation for the Wigner semicircle law (Lemma \ref{lm:stabilitymsc}) which is defined as
$$
m_{sc}(z)=\int_{\R} \frac{1}{x-z}\d\mu_{sc}(x),\qquad \mu_{sc}(x)=\frac{1}{2\pi} \sqrt{4-x^2}\mathbbm 1_{-2\leq x\leq 2}.
$$
 (One can also prove stability for equation \eqref{eq:mpquadratic'} directly; see \cite{BKYY}.) 

\be{lm}[cf.\ Lemma 7.6 in \cite{EY}]\label{lm:stabilitymsc}
Let $z=E+i\eta$ with $|E|\leq 20$, $0<\eta\leq 10$ and $\kappa=||E|-2|$. Suppose that $m$ satisfies
$$
\l|m+\frac{1}{z+m}\r|\leq \de
$$
for some $\de\leq 1$. Then
$$
\min\l\{|m-m_{sc}(z)|,\l|m-\frac{1}{m_{sc}(z)}\r|\r\}\leq \frac{C\de}{\sqrt{\kappa+\eta+\de}}
$$
\e{lm}

\be{proof}[Proof of Lemma \ref{lm:stability}]
We use the substitution
\beq\label{eq:sqrtz}
\tilde z= \sqrt{z},\qquad \tilde m=\sqrt{z}m,\qquad \tilde \delta=\de|\sqrt{z}|,
\eeq
where $\sqrt{z}$ is the branch of the square root defined by $\sqrt{re^{i\theta}}=\sqrt{r} e^{i\theta/2}$ for all $\theta\in (-\pi,\pi)$. Multiplying \eqref{eq:stabilityass} by $\sqrt{z}$ and substituting gives
$$
\l|\tilde m+\frac{1}{\tilde z+\tilde m}\r|\leq \tilde\de.
$$
We may verify that $z\in\curly{D}$ ensures that $\tilde{z}=\sqrt{z}=\tilde E+i\tilde \eta$ has $\tilde E,\tilde\eta $ satisfying the assumptions in Lemma \ref{lm:stabilitymsc}. Moreover, we note
that for $z\in\curly{D}$, the identity $|\sqrt{z}|=\sqrt{|z|}$, the assumption $\eta\leq c_\kappa$, the fact that we can assume without loss of generaliy that $c_\kappa\leq \kappa$, and the inequality $(1-x)^{1/4}\leq 1-x/4$ imply
$$
\begin{aligned}
&\tilde \kappa =||\tilde E|-2|=2-\tilde E
\geq 2-|\sqrt{z}|
\geq 2-((4-\kappa)^2+\kappa)^{1/4}\\
&\geq 2\l(1-\l(1-\frac{3}{8}\kappa\r)^{1/4}\r)\geq \frac{3}{16}\kappa.
\end{aligned}.
$$
Now we apply Lemma \ref{lm:stabilitymsc} and divide by $|\sqrt{z}|$ afterwards to conclude
\beq\label{eq:intermediate}
\min\l\{\l|m-\frac{m_{sc}(\sqrt{z})}{\sqrt{z}}\r|,\l|m-\frac{1}{\sqrt{z} m_{sc}(\sqrt{z})}\r|\r\}\leq \frac{C\de}{\sqrt{\tilde\kappa+\eta+\de}}
\leq \frac{C}{\sqrt{\kappa}}.
\eeq

To prove the claim, it remains to show that $\frac{m_{sc}(\sqrt{z})}{\sqrt{z}}=m_{\mathrm{\mathrm{MP}}}(z)$. This can be verified directly from the definition of these Stieltjes transforms via the substitution $x=\sqrt{y}$, the fact that $\mu_{sc}(x)=\mu_{sc}(-x)$, and a partial fraction decomposition. Lemma \ref{lm:stability} then follows from \eqref{eq:intermediate}.
\e{proof}


\subsection{Order-one bounds on $m_{\mathrm{\mathrm{MP}}}(z)$}
For later use, we recall the following well-known bounds on $m_{\mathrm{\mathrm{MP}}}(z)$. 

\be{lm}
\label{lm:mMPbounds}
There exist constants $C_\kappa,C_{\kappa}'>1$ so that
\beq\begin{aligned}
\frac{1}{C_\kappa} \leq |m_{\mathrm{\mathrm{MP}}}(z)|\leq C_{\kappa},\qquad
\frac{1}{C_\kappa'}  \leq \Im[m_{\mathrm{\mathrm{MP}}}(z)]\leq C_{\kappa}',\qquad \forall z\in\curly{D}.
\end{aligned}
\eeq
\e{lm}

\be{proof}
These bounds are straightforward consequences of the explicit formula
$$
m_{\mathrm{\mathrm{MP}}}(z)=\frac{-z+i \sqrt{z(4-z)}}{2z}
$$
and our definition of the domain $\curly{D}$. See Lemma 3.3 in \cite{BKYY} for more details.
\e{proof}

These bounds allow us to simplify the analysis in the spectral bulk, our main area of interest, by noting that the second term in the minimum in Lemma \ref{lm:stability} is always of order $1$, i.e., large. 

\be{cor}\label{cor:imlb}
We have
\beq\label{eq:otherterm}
\l|m_N(z) - \frac{1}{zm_{\mathrm{\mathrm{MP}}}(z)}\r|
\geq C_{\kappa}.
\eeq
for all $z\in\curly{D}$ provided that the constant $c_\kappa$ in the definition of $\curly{D}$ is sufficiently small.
\e{cor}

\be{proof}
By Lemma \ref{lm:mMPbounds} and the definition of $\curly{D}$, 
$$
\Im[zm_{\mathrm{\mathrm{MP}}}(z)]=\Re[z]\Im[m_{\mathrm{\mathrm{MP}}}(z)]+\Im[z]\Re[m_{\mathrm{\mathrm{MP}}}(z)]
\geq \frac{\kappa}{C_\kappa'}-c_\kappa C_\kappa
$$
and this equals a positive constant $C_\kappa$ for sufficiently small $c_\kappa$. Thus
$$
\l|m_N(z) - \frac{1}{zm_{\mathrm{\mathrm{MP}}}(z)}\r|\geq
\Im\l[m_N(z) - \frac{1}{zm_{\mathrm{\mathrm{MP}}}(z)}\r]
\geq \frac{\Im[zm_{\mathrm{\mathrm{MP}}}(z)]}{|zm_{\mathrm{\mathrm{MP}}}(z)|^2}
\geq C_{\kappa,\eps}
$$
where we applied Lemma \ref{lm:mMPbounds} again to the denominator in the last step. This proves Corollary \ref{cor:imlb}.
\e{proof}

\subsection{Initial scale estimate}
At this point, we can establish the main conclusion for sufficiently large $\eta$. We recall the definition of $\theta_0$ from Theorem \ref{thm:main1}. It is convenient to introduce 
\beq\label{eq:beta0defn}
\beta_0:=4\theta_0+1
\eeq

\be{prop}[Initial scale estimate]\label{prop:step4}
Let $z=E+i\eta\in \curly{D}$ with $\eta= N^{-\theta/4}$ and $0<\theta<\theta_0$. Then
\begin{equation}
    \mathbb P\l(|m_N(z) - m_{\mathrm{\mathrm{MP}}}(z)| > \frac{N^{-\theta_0}}{\eta}\r)\leq 
    N^{1-\beta_0}.
\end{equation}
\e{prop}


The proof uses the following lemma wherein we use the same branch cut for $\sqrt{z}$ as in \eqref{eq:sqrtz}. 
\be{lm}\label{lm:imsign}
Let $\Im[z]>0$. Then $\Im[\sqrt{z}m_N(z)]>0$.
\e{lm}

This lemma is motivated by the identity $\sqrt{z}m_{\mathrm{\mathrm{MP}}}(z)=m_{sc}(\sqrt{z})$. The crucial observation is that $m_{N}$ is the Stieltjes transform of a measure supported on $\R_+$.  

\be{proof}[Proof of Lemma \ref{lm:imsign}]
We use a partial fraction decomposition to write
$$
\sqrt{z}m_N(z)
=\sqrt{z}\frac{1}{N}\sum_{j=1}^N \frac{1}{\sigma_j^2-z}
=\frac{1}{2N}\sum_{j=1}^N \l(\frac{1}{\sigma_j-\sqrt{z}}-\frac{1}{\sigma_j+\sqrt{z}}\r)
$$ 
and the latter expression has positive imaginary part whenever $\Im[z]>0$.
\e{proof}

We are now ready to give the

\be{proof}[Proof of Proposition \ref{prop:step4}]
Let $z\in\curly{D}$. From \eqref{eq:step2} and elementary estimates, we obtain
\beq\label{eq:step4}
\begin{aligned}
\l|m_N(z)+\frac{1}{z+zm_N(z)}\r|
\leq& \frac{1}{N} \sum_{i=1}^N \l|\frac{1}{z+zm_N(z)+\curly{E}_i(z)}-\frac{1}{z+zm_N(z)} \r|\\
\leq& 2\frac{1}{|z+zm_N(z)|^2} \max_{1\leq i\leq N}|\curly{E}_i(z)|
=:\de
\end{aligned}
\eeq
provided that
\beq\label{eq:Omsmall}
\max_{1\leq i\leq N}|\curly{E}_i(z)|\leq \frac{1}{2}|z+zm_N(z)|.
\eeq

To ensure \eqref{eq:Omsmall}, we first estimate the right-hand side using Lemma \ref{lm:imsign} 
\beq\label{eq:initiallb}
\frac{1}{2}|z+zm_N(z)|
\geq \frac{1}{2}|\sqrt{z}|\Im[\sqrt{z}+\sqrt{z}m_N(z)]
\geq \frac{\sqrt{\kappa}}{2} \Im[\sqrt{z}]
\eeq
In polar coordinates, since $z\in\curly{D}$,
\beq\label{eq:initiallb'}
\Im[\sqrt{z}]=(E^2+\eta^2)^{1/4}\sin\arctan\l(\frac{\eta}{2E}\r)
\geq C\sqrt{\kappa} \eta.
\eeq
Hence, the condition\eqref{eq:Omsmall} is implied by the stronger condition
\beq\label{eq:Omsmall'}
\max_{1\leq i\leq N}|\curly{E}_i(z)|\leq \frac{\kappa}{8\pi}\eta.
\eeq
To verify \eqref{eq:Omsmall'}, we estimate each $|\curly{E}_i(z)|\leq |m_N(z)-m_N^{(i)}(z)|+|F_i(z)|$ via Lemma \ref{lm:interlace} and Corollary \ref{cor:ld} with $\beta=\beta_0$ given by \eqref{eq:beta0defn}. From these and a union bound, we obtain that
\beq\label{eq:maxOm}
\max_{1\leq i\leq N}|\curly{E}_i(z)|\leq \frac{C}{N\eta}+\frac{N^{-\theta_0}}{\eta}
\leq C\frac{N^{-\theta_0}}{\eta}
\leq \frac{\kappa}{8\pi}\eta
\eeq
holds except on a set of probability $\leq N^{1-\beta_0}$. In the last step we used that $\eta=N^{-\theta/4}$ with $0<\theta<\theta_0$. The upshot of these considerations is that \eqref{eq:step4} holds except on a set of probability $\leq N^{1-\beta_0}$. 

Let us therefore assume that \eqref{eq:step4} holds. We apply Lemma \ref{lm:stability}, the bounds in \eqref{eq:initiallb}, \eqref{eq:initiallb'}, and \eqref{eq:maxOm} to find
$$
\min\left\{|m_N(z) - m_{\mathrm{\mathrm{MP}}}(z)|,  \l|m_N(z) - \frac{1}{zm_{\mathrm{\mathrm{MP}}}(z)}\r|\right\} 
\leq C \frac{\delta}{\kappa}
\leq C \frac{N^{-\theta_0}}{\eta^3}
\leq C\eta
$$
where the last step uses $\eta=N^{-\theta/4}$ with $0<\theta<\theta_0$. In particular, we see that the right-hand side vanishes as $N\to\infty$, while the second term on the left-hand side is bounded below by a positive constant uniform in $N$ by Corollary \ref{cor:imlb}. Hence,
\beq\label{eq:mfirstbound}
|m_N(z) - m_{\mathrm{\mathrm{MP}}}(z)|\leq C \frac{\delta}{\kappa}\leq C \eta
\eeq
 With this bound in hand, we can improve the a priori lower bounds \eqref{eq:initiallb} and \eqref{eq:initiallb'} to order-$1$ constants, i.e.,
$$
\begin{aligned}
|z+zm_N(z)|\geq& |z+ z m_{\mathrm{MP}}(z)|- |z m_{\mathrm{MP}}(z) - z m_N(z)|\geq \Im[zm_{\mathrm{\mathrm{MP}}}(z)]-C\eta\\
 \geq& \Re[z]\Im[m_{\mathrm{\mathrm{MP}}}(z)]-C\eta\geq C
\end{aligned}
$$
where the last step uses Lemma \ref{lm:mMPbounds}. Together with \eqref{eq:maxOm}, this implies that 
$$
\de=2\frac{\max_{1\leq i\leq N}|\curly{E}_i(z)|}{|z+zm_N(z)|^2} 
\leq C\max_{i}|\curly{E}_i(z)|\leq C \frac{N^{-\theta_0}}{\eta}.
$$
 Replacing the last estimate in \eqref{eq:mfirstbound} with this improved bound on $\de$, we conclude that for all $z\in\curly{D}$,
$$
|m_N(z) - m_{\mathrm{\mathrm{MP}}}(z)|\leq C\de\leq  C\frac{N^{-\theta_0}}{\eta}
$$
holds except on a set of probability $\leq N^{1-\beta_0}$. This proves Proposition \ref{prop:initial}.
\e{proof}

\section{Step 5: Bootstrap argument and conclusion}
The allowed values of $\eta$ can be improved from Proposition \ref{prop:initial} down to the scale $\eta=N^{-\theta_0}$ that is seen in the main result. This uses a by now standard bootstrap argument based on Lipschitz continuity. Here we use a straightforward modification of the standard argument that leads to significantly better constraints on the main parameters $\theta$ and $d$. 

The simple observation is that the derivative $|m'_N(z)|\leq \eta^{-2}\leq N^{2\theta}$ for our purposes. Since the bootstrap argument relies on the mean-value theorem in the form
\beq\label{eq:mvt}
|m_N(z_1)-m_N(z_2)|\leq |z_1-z_2|N^{2\theta},\qquad \text{for } z_1,z_2\in \curly{D},
\eeq 
this trivial refinement allows us to use a coarser lattice spacing than the usual $N^{-4}$ \cite{EY}. Consequently, we only require a union bound for the probability over a relatively small collection of events and can get by with the relatively small concentration exponent $\beta_0$ from \eqref{eq:beta0defn}.

The following proposition summarizes a single step in the bootstrap argument. For $c>0$ and $z\in\curly{D}$, we define the event
$$
\Om_{c}(z):=\l\{|m_N(z)-m_{\mathrm{\mathrm{MP}}}(z)|\leq c\frac{N^{-\theta_0}}{\Im[z]}\r\}.
$$

\be{prop}[Bootstrap argument]\label{prop:bootstrap}
For every sufficiently large constant $c_\kappa'>0$ the following holds. Let $z_1=E+i\eta_1\in \curly{D}$ and let $\Om$ be an event satisfying 
$$
\Om\subset \Om_{c_\kappa'}(z_1).
$$
Let $z_2=E+i\eta_2$ with $|\eta_1-\eta_2|\leq N^{-s}$ for some $s>\theta+\theta_0$. Then
$$
\mathbb P\l(\Om\setminus \Om_{c_\kappa'}(z_2)\r)\leq N^{1-\beta_0}
$$
holds for all sufficiently large $N$.
\e{prop}

\be{proof}
Suppose that $\Om$ occurs. Then we have $|m_N(z)-m_{\mathrm{\mathrm{MP}}}(z)|\leq c_\kappa'\frac{N^{-\theta_0}}{\Im[z_1]}$. We aim to use Lemma \ref{lm:stability} with $z=z_2$ to establish $\Om_{c_\kappa'}(z_2)$ except on a set of small probability. To use Lemma \ref{lm:stability} effectively, we first need to control $|z_2+z_2m_N(z_2)|^{-1}$. The triangle inequality, the mean-value theorem \eqref{eq:mvt}, and its analog for $m_{\mathrm{\mathrm{MP}}}$ imply
$$
\begin{aligned}
&|m_N(z_2)-m_{\mathrm{\mathrm{MP}}}(z_2)|\\
&\leq |m_N(z_2)-m_N(z_1)|+|m_N(z_1)-m_{\mathrm{\mathrm{MP}}}(z_1)|+|m_{\mathrm{\mathrm{MP}}}(z_1)-m_{\mathrm{\mathrm{MP}}}(z_2)|\\
&\leq 2N^{2\theta-s}+c_\kappa' N^{\theta-\theta_0}.
\end{aligned}
$$
Since $\theta<\theta_0$ and $s>2\theta_0$, the right-hand side is $o(1)$ as $N\to\infty$. Hence, by the triangle inequality and Lemma \ref{lm:mMPbounds},
\beq\label{eq:eff}
|z_2+z_2m_N(z_2)|\geq |z_2||1+m_{\mathrm{\mathrm{MP}}}(z_2)|-|z_2||m_N(z_2)-m_{\mathrm{\mathrm{MP}}}(z_2)|\geq 2C_\kappa-o(1)\geq C_\kappa
\eeq
for all sufficiently large $N$. This is the required control on $|z_2+z_2m_N(z_2)|^{-1}$.

We apply Corollary \ref{cor:ld} with $\beta=\beta_0$ from \eqref{eq:beta0defn} and a union bound to ensure that $\max_{1\leq i\leq N}|\curly{E}_i(z)|\leq |z+zm_N(z)|$ and in fact
\beq\label{eq:debound}
\de=2\frac{\max_{1\leq i\leq N}|\curly{E}_i(z_2)|}{|z_2+z_2m_N(z_2)|^2} 
\leq C\max_{i}|\curly{E}_i(z_2)|\leq C \frac{N^{-\theta_0}}{\eta}
\eeq
holds except on a set of probability $\leq N^{1-\beta_0}$. From now on, we assume that \eqref{eq:debound} holds. By Lemma \ref{lm:stability} and \eqref{eq:debound}, 
$$
\min\left\{|m_N(z_2) - m_{\mathrm{\mathrm{MP}}}(z_2)|,  \l|m_N(z_2) - \frac{1}{z_2m_{\mathrm{\mathrm{MP}}}(z_2)}\r|\right\} 
\leq C \de\leq  C \frac{N^{-\theta_0}}{\eta}.
$$
Since the right-hand side vanishes as $N\to\infty$, Corollary \ref{cor:imlb} implies that
$$
|m_N(z_2) - m_{\mathrm{\mathrm{MP}}}(z_2)|\leq C \frac{N^{-\theta_0}}{\eta}.
$$
We have shown that $\Om_{c_\kappa'}(z_2)$ occurs except on a set of probability $\leq N^{1-\beta_0}$. This proves Proposition \ref{prop:bootstrap}.
\e{proof}

\subsection{Conclusion}
We are now ready to give the 

\be{proof}[Proof of Theorem \ref{thm:main1}]
We set $s=\frac{\theta}{2}+\frac{3\theta_0}{2}>\theta+\theta_0$ and we discretize the domain $\curly{D}$ into the lattice
\beq\label{eq:tildeDdefn}
\tilde{\curly{D}}:=N^{-s}(\Z+i(\Z+N^{-\theta/4}))\cap \curly{D}
\eeq
where $\Z+N^{-\theta/4}$ are the integers shifted by $N^{-\theta/4}$. We first apply Proposition \ref{prop:step4} to every $z$ that lies in the intersection of $\tilde{\curly{D}}$ and the line $\curly{L}=\{z\in\C\,:\,\Im[z]=N^{-\theta/4}\}$.  By a union bound over the order $N^s$ many exceptional events, this implies that 
\beq\label{eq:mdiffclaim1}
|m_N(z)-m_{\mathrm{\mathrm{MP}}}(z)|\leq c_\kappa' \frac{N^{-\theta_0}}{\Im[z]}
\eeq
holds for all $z\in \tilde{D}\cap \curly{L}$ except on a set of probability $\leq C N^{s+1-\beta_0}$. Then, we apply Proposition \ref{prop:bootstrap} to cover all other possible imaginary values in $\curly{D}$. Altogether, this requires a total order of $N^{2s}$-many union bounds. The upshot is that \eqref{eq:mdiffclaim1} holds for all $z\in \tilde{\curly{D}}$ with probability $\leq C N^{2s+1-\beta_0}$. We note that $2s+1-\beta_0=\theta-\theta_0<0$ for our choice of parameters $\beta_0=4\theta_0+1$ and $s=\frac{\theta}{2}+\frac{3\theta_0}{2}$. Hence, we have shown that \eqref{eq:mdiffclaim1} holds for all $z\in \tilde{\curly{D}}$ with high probability in the sense of Definition \ref{defn:whp}.

Finally, by Lipschitz continuity of $m_N$ and $m_{\mathrm{\mathrm{MP}}}$ in the form of \eqref{eq:mvt} and the triangle inequality, we conclude that for all $z\in \curly{D}$,
$$
|m_N(z)-m_{\mathrm{\mathrm{MP}}}(z)|
\leq c_\kappa' \frac{N^{-\theta_0}}{\Im[z]}+2N^{2\theta-s}
\leq 2c_\kappa' \frac{N^{-\theta_0}}{\Im[z]},
$$
for all sufficiently large $N$. The last step uses that $s>\theta+\theta_0$. Note that the mean-value theorem applies deterministically and hence without further loss in probability. This proves Theorem \ref{thm:main1}. 
 \e{proof}
 
 \section*{Acknowledgments}
 
 We are immensely grateful to Horng-Tzer Yau for pointing us to the useful operator identity \eqref{eq:B*B}. 
 
\begin{appendix}
\section{Proof of Corollary \ref{cor:rigidity}}
By Theorem \ref{thm:main1}, it suffices to prove that \eqref{eq:main} implies \eqref{eq:corcounting}.
This is a standard argument based on the Helffer-Sj\"ostrand formula with only minor modifications to adapt to the scale $N^{-\theta}$ with $\theta<\theta_0$. More precisely, considering Lemma 11.2 in \cite{EY}, we only need to observe that the constants $U_1$ and $U_2$ are less than $N^{-\theta}$. For $U_1$ this holds by definition. The computation for $U_2$ was performed in Equation 11.36 in \cite{EY} using the monotonicity of $\sqrt{z} \Im[m_N(z)]$ and $\sqrt{z} \Im[m_{\mathrm{\mathrm{MP}}}(z)]$ in $z$ and the errors are smaller than $N^{-\theta}$. The details are left to the interested reader who may also find it helpful to refer to the proof of Lemma 11.3 in \cite{EY}. 
\qed

\section{Proof of Theorem \ref{thm:deloc}}
 Fix $1\leq \al,i\leq N$ and let $z=E+i\eta\in\tilde{\curly{D}}$ with the lattice $\tilde{\curly{D}}$ defined in \eqref{eq:tildeDdefn} and $\eta=CN^{-\theta}$ for $\theta\in (0,\theta_0)$. Notice that 
\beq\label{eq:evectorinitial}
 |u_\al(i)|^2
\leq \frac{\eta^2}{\eta^2+(\sig_j^2-E)^2} |u_\al(i)|^2
=
 \eta \Im \l(\frac{1}{X_NX_N^*-z}\r)_{ii}
\eeq
The resolvent can be expressed via the Schur complement formula, cf.\  \eqref{eq:Schurcomplement'}, as
$$
\Im \l(\frac{1}{X_NX_N^*-z}\r)_{ii}
=\frac{-1}{z+z m_N^{(i)}(z)+F_i}
$$
with $F_i$ defined in \eqref{eq:Fidefn}.

By \eqref{eq:debound}, Lemma \ref{lm:interlace} and Theorem \ref{thm:main1}, we know that 
\beq\label{eq:maxappears}
\max\l\{\max_{1\leq i\leq N}|F_i|,|m_N^{(i)}(z)-m_{\mathrm{\mathrm{MP}}}(z)|\r\}\leq C\frac{N^{-\theta_0}}{\eta}
\eeq
holds with high probability. By these facts and the fact that $m_{\mathrm{\mathrm{MP}}}(z)$ is an order-$1$ quantity in the sense established by Lemma \ref{lm:mMPbounds}, we find that
$$
\l|\frac{-1}{z+z m_N^{(i)}(z)+F_i}\r|\leq C
$$
Hence, \eqref{eq:evectorinitial} implies
$$
|u_\al(i)|^2\leq C\eta=C N^{-\theta}.
$$
In view of \eqref{eq:maxappears}, this estimate is uniform in $i$ and so we can take the maximum over $i\in \{1,\ldots,N\}$. This  proves Theorem \ref{thm:deloc}.
\qed
\mbox{}
\end{appendix}

\mbox{}

\end{document}